\documentclass[12pt]{article}
\usepackage{amsfonts, amssymb,version}
\def\bb{{\bf b}}
\def\C{\mathbb C}

\def\dna{\downarrow}

\def\EE{{\bf E}} 
 
\def\FF{\mathbb F}

\def\LL{{\bf L}}

\def\M{{\cal M}} 
\def\MM{{\bf M}} 

\def\N{{\cal N}} 
 
\def\O{{\cal O}}

\def\PP{{\bf P}}

\def\QQ{{\bf Q}}
\def\R{\mathbb R}

\def\t{\tau}
\def\T{{\cal T}}
\def\U{{\cal U}}

\def\Z{\mathbb Z}

\def\deg{\mathop{\rm deg}\nolimits} 
 
\def\id{\mathop{\rm id}\nolimits} 
\def\im{\mathop{\rm im}\nolimits} 
\def\ker{\mathop{\rm ker}\nolimits}

\def\rank{\mathop{\rm rank}\nolimits}

\def\sq{\mathop{\rm Sq}\nolimits}
\let\hra\hookrightarrow

\let\lra\longrightarrow

\let\ov\overline
\def\ovnu{\overline{\nu}}

\let\wh\widehat
\let\pa\partial

\newtheorem{theorem}{Theorem}[section]
\newtheorem{proposition}[theorem]{Proposition}
\newtheorem{lemma}[theorem]{Lemma}

\newtheorem{corollary}[theorem]{Corollary}

\def\rem{\refstepcounter{theorem}\paragraph{Remark \thetheorem}}

\def\proof{\paragraph{Proof}}
\def\example{\refstepcounter{theorem}\paragraph{Example \thetheorem}}
\def\defin{\refstepcounter{theorem}\paragraph{Definition \thetheorem}}

\makeatletter \def\l@section{\@dottedtocline{1}{0em}{1.2em}} \makeatother

\begin{document}

\title{Topology of Quadric Bundles}
\author{Yogish I. Holla and Nitin Nitsure}
\date{}

\maketitle

\tableofcontents

%
\section{Introduction}

The complex general orthogonal group $GO(r)$ 
by definition consists of all $r \times r$  
complex matrices $g$ such that ${^t}g g = c I$ where $c$ is an invertible
scalar. (For any integer $r\ge 1$, we will write $GO(r,\C)$, $GL(r,\C)$, 
$O(r,\C)$ simply as $GO(r)$, $GL(r)$ and $O(r)$, omitting the `$\C$\,'.)
A principal $GO(r)$-bundle on any base space $X$ is then equivalent to 
a {\bf non-degenerate triple} $\T= (E,L,b)$ consisting of a rank $r$ complex
vector bundle $E$, a complex line bundle $L$ on $X$ and a  
symmetric bilinear form $b : E\otimes E \to L$ on $E$ 
taking values in $L$, such that $b$ is non-degenerate on
each fiber. The vector bundle $E$ corresponds to the defining representation
$GO(r)\hra GL(r)$ of $GO(r)$, the line bundle  corresponds to the 
character $\sigma :GO(r) \to \C\,^*$ defined by the equation
${^t}g g = \sigma(g) I$, and $b$ is induced by the standard
bilinear form $\sum x_iy_i$ on $\C\,^r$. 
The reverse correspondence is by a 
Gram-Schmidt process applied locally (\'etale-locally in the algebraic
category). We define a {\bf non-degenerate quadric bundle} $Q$ of rank
$r$ on a base $X$ to be an equivalence class $[\T]$ of non-degenerate rank 
$r$ triples $\T$, where for any line bundle $K$ on $X$, 
we regard the triples $\T$ and $\T\otimes K = (E\otimes K,L\otimes K^2,b
\otimes 1_K)$ to be equivalent.

The central result (called the Main Theorem) is proved
for cohomology with coefficients in an arbitrary ring. The computational
results and examples involving cohomological classes are formulated 
and proved for the case of $\FF_2 = \Z/(2)$ coefficients.
We always work with singular cohomology. 

A cohomological {\bf quadric invariant} $\alpha$ of 
non-degenerate quadric bundles of rank $r$ associates to any
such quadric bundle $Q$ on any base $X$ an element $\alpha(Q)\in H^*(X)$,
with $f^*(\alpha(Q)) = \alpha(f^*Q)\in H^*(X')$ 
under any pullback $f : X'\to X$. Let $H^*(BGO(r))$ 
denote the singular cohomology ring with coefficients
$\FF_2$ of the classifying space $BGO(r)$ of $GO(r)$. 
It can be seen that the cohomological invariants of non-degenerate quadric
bundles of rank $r$ bijectively 
correspond to elements of the so called {\bf primitive subring}
$PH^*(BGO(r))\subset H^*(BGO(r))$,
which consists of all elements $x \in H^*(BGO(r))$ such that 
$B(\mu)^*(x)=1\otimes x$.
$B(\mu)^*: H^*(BGO(r))\to H^*(B\C^*)\otimes H^*(BGO(r))$
induced by the multiplication map $\mu : \C^*\times GO(r) \to GO(r)$, 
which corresponds to the operation $\T \mapsto \T \otimes K$ on triples.

More generally, we consider triples $\T= (E,L,b)$ which are 
not assumed to be non-degenerate, but for which for each $x\in X$, 
the restriction $b_x : E_x\otimes E_x \to L_x$ to the fiber 
over $x$ is non-zero.
Again defining an equivalence relation $\T \sim \T \otimes K$, 
we call an equivalence class $Q=[\T]$ as a {\bf quadric bundle} on $X$.
We say that the triple $\T$, or the corresponding quadric bundle $Q$, 
is {\bf minimally degenerate} if the rank of $b_x$ is constant equal to
$r-1$ for all $x\in X$, where $r= \rank(E)$. Note that such a $Q$ 
naturally defines a non-degenerate quadric bundle $\ov{Q}$ of rank $r-1$
on $X$, where $r=\rank(Q)$. If a triple $(E,L,b)$ represents $Q$, 
then $\ov{Q}$ is represented by the triple $(\ov{E}, L, \ov{b})$ where
$\ov{E} = E/\ker(b)$ and $\ov{b}$ is induced by $b$.

Given a Hausdorff topological space $X$ of the homotopy type of a CW-complex,
we say that a closed subspace $Y\subset X$ of the homotopy type of a 
CW-complex is a {\bf topological divisor} if there exists an
open neighbourhood $U$ of $Y$ in $X$, such that $U$ is homeomorphic to
the total space of a real rank $2$ vector bundle $N$ on $Y$, under a 
homeomorphism which is identity on $Y$
(where we regard $Y$ as the zero section of $N$). 
An important example of this is when
$X$ is a complex algebraic variety, and $Y\subset X$ is a 
smooth Weil divisor, which does not meet the singular locus of $X$.
Given such a pair $(X,Y)$, we will
consider quadric bundles $Q=[\T]$ on $X$ which are non-degenerate
on $X-Y$, and are minimally degenerate on $Y$. We will say that such a 
$Q$ is a {\bf mildly degenerating} quadric bundle on $(X,Y)$, of
{\bf generic rank} $r=\rank(E)$. 

This paper is devoted to the following three questions.

{\bf (1) } Let the coefficients be $\FF_2$.
What is the cohomology ring $H^*(BGO(r))$, in terms of generators and
relations? Under the canonical inclusions $O(r) \hra GO(r)$ and
$GO(r) \hra GL(r)$,
what are the induced ring homomorphisms $H^*(BGL(r)) \to H^*(BGO(r))$
and $H^*(BGO(r)) \to H^*(BO(r))$? 

{\bf (2) } Let the coefficients be $\FF_2$.
What are the invariants of non-degenerate quadric bundles? Describe 
concretely the primitive subring $PH^*(BGO(r))\subset H^*(BGO(r))$.

{\bf (3) }Let the coefficients be in an arbitrary ring $R$. 
Consider a space $X$ and together with a  topological divisor
$Y\subset X$. For any mildly degenerating quadric bundle $Q$ on $(X,Y)$
of generic rank $r$,
how are the quadric invariants of the restriction $Q_{X-Y}$
related to the quadric invariants
of the rank $r-1$ non-degenerate quadric bundle $\ov{Q_Y}$ 
corresponding to the restriction $Q_Y$ ? 

We now sketch the answers.

When $r = 2n+1$ is odd, the first two questions are easily answered,
using the group isomorphism $\C^* \times SO(2n+1) \to GO(2n+1)$
given by $(a,g) \mapsto ag$, and the known facts about 
Stiefel-Whitney and Chern classes. The details are given in
Section 2.

When $r=2n$ is even, the ring $H^*(BGO(2n))$ has been described by means of
finitely many generators $\lambda,\, a_{2i-1}, \, d_T,\, b_{4i}$ 
and relations between these generators (see [H-N]).
The ring homomorphism $H^*(BGO(2n)) \to H^*(BO(2n))$ is also determined in
[H-N]. We determine the ring homomorphisms $H^*(BGL(2n)) \to H^*(BGO(2n))$
in Section 4 (Proposition \ref{odd Chern classes}), which completes the
answer to question (1) in even ranks. 

The determination of the map $H^*(BGL(2n)) \to H^*(BGO(2n))$ 
enables us to write the ring homomorphism 
$B(\mu)^*: H^*(BGO(2n))\to H^*(B\C^*)\otimes H^*(BGO(2n))$ in terms of the
generators $\lambda,\, a_{2i-1}, \, d_T,\, b_{4i}$.
Consider the similar ring homomorphism 
$B(\mu)^*: H^*(BO(2n))\to H^*(B(\Z/(2)))\otimes H^*(BO(2n))$ induced
by the multiplication map $\mu : \Z/(2) \times O(2n) \to O(2n)$ where
$\Z/(2) = \{ \pm 1\}$. 
When the rank $r=4m+2$ is congruent to $2$ mod $4$, 
and also in the case $r=4$, H. Toda has determined 
the corresponding primitive ring $PH^*(BO(r))\subset H^*(BO(r))$
(see Toda [T]).
This covers (one more than!) half the even cases. 
Using these results of Toda, we are able to give a finite set of generators 
for the primitive rings $PH^*(BGO(4m+2))\subset H^*(BGO(4m+2))$ and also
for $PH^*(BGO(4))$. For the remaining ranks $8, 12, 16, \ldots$, we have not
fully answered the question (2). However, even in these cases, 
our concrete description of the map 
$B(\mu)^*$ makes it possible to mechanically check 
whether any given polynomial 
$f(\lambda,\, a_{2i-1}, \, d_T,\, b_{4i})\in H^*(BGO(2n))$, 
is in $PH^*(BGO(2n))$.

Now we come to the question (3), the answer 
to which is our main result. 
First observe that if $Q = [E,L,b]$ is a 
minimally degenerate quadric bundle on
any base $Y$ (that is, $\rank(b_y) = \rank(E) -1$ at all $y\in Y$)
then (up to isomorphism) we have a {\bf canonical triple} 
$\T^Q$ which represents the non-degenerate rank $r-1$ quadric bundle
$\ov{Q}$, defined as follows. Let $(E,L,b)$ be any 
arbitrary triple with $Q = [E,L,b]$, so that 
$\ker(b : E \to L\otimes E^*)$ is a line subbundle of $E$.
Let $\ov{E}$ be the quotient $E/\ker(b)$, and let
$\ov{b} : \ov{E}\otimes\ov{E}\to L$ be induced by
$b$. Note that $(\ov{E},L,\ov{b})$ is a non-degenerate rank $r-1$ triple,
where $r=\rank(E)$, which defines $\ov{Q}$. 
We define the rank $r-1$ triple $\T^Q$ as 
$$\T^Q = (\ov{E}, L, \ov{b}) \otimes (\ker(b))^{-1}$$
This is well-defined, independent of the initial choice of a 
representative $(E,L,b)$ for $Q$. 
In particular, if $Y\subset X$ is a 
topological divisor and $Q$ is 
a mildly degenerating quadric bundle on $(X,Y)$ of generic rank $r$,
then we get a canonically defined 
principal $GO(r-1)$-bundle on $Y$ corresponding to the triple $\T^{Q_Y}$. 

Next, observe that it follows from the definition of a topological
divisor $Y\subset X$ that we have a Gysin boundary map 
$$\pa : H^*(X-Y) \to H^{*-1}(Y)$$ 
which is an additive map, homogeneous of degree $-1$.  
Our main theorem says that for any quadric invariant 
$\alpha(Q_{X-Y}) \in H^*(X-Y)$ of the 
non-degenerate rank $r$ quadric bundle 
$Q_{X-Y}$, the Gysin boundary $\pa(\alpha(Q_{X-Y})) \in H^*(Y)$
is expressible in terms of the $GO(r-1)$-characteristic classes of the 
triple $\T^{Q_Y}$ as follows. First we need some preliminaries.

If $Y\subset X$ is a topological divisor, and $K$ is a line bundle on
$X$ together with a section $s \in \Gamma(X-Y,K)$ which is everywhere 
non-vanishing on $X-Y$, we can define its 
{\bf absolute topological vanishing multiplicity} 
$|\nu|_Y(s)\in H^0(Y;\Z)$, and its {\bf topological vanishing parity}
$\ovnu_Y(s) \in H^0(Y;\Z)$ 
(see Section 6). 
In the algebraic category, this is the usual vanishing multiplicity
$\nu_Y(s)$ and its parity, where even is $0$ and 
odd is $1$.
Given a quadric bundle $Q=[E,L,b]$ on $X$
which is non-degenerate outside a topological divisor $Y\subset X$,
the discriminant 
$\det(b) \in \Gamma(X, L^{rank(E)} \otimes \det(E)^{-2})$, 
and its topological vanishing multiplicity 
$\nu_Y(\det(b))$ along $Y$, are independent of the choice of 
a triple $(E,L,b)$ that represents the 
given quadric bundle $Q$, which allows us to denote these
simply as $\det(Q) \in \Gamma(X, L^r\otimes \det(E)^{-2})$ and as 
$\nu_Y(\det(Q))$.  

Let $B(v)^* : H^*(BGO(r)) \to  H^{*-1}(BO(r-1))$ be the ring
homomorphism induced by the group homomorphism $v: O(r-1) \to GO(r)$
defined by $g\mapsto \left(\begin{array}{cc} 1& \\ & g \end{array}\right)$.
Let $(\EE,\LL,\bb)$ be the universal triple over $BGO(r-1)$.
The complement $\LL_o$ of the zero section of $\LL$ can be taken to be
$BO(r-1)$ (see Section 2), so there is a Gysin boundary map
$d : H^*(BO(r-1)) \to H^*(BGO(r-1))$. 
Finally, let $\delta : PH^*(BGO(r)) \to H^{*-1}(BGO(r-1))$ be the
composite linear map 
$$PH^*(BGO(r)) \hra H^*(BGO(r)) \stackrel{B(v)^*}{\to} H^*(BO(r-1))
\stackrel{d}{\to} H^{*-1}(BGO(r-1))$$
Here the cohomologies have coefficients in some fixed ring $R$. 
With the above notations, we can now state the following.

\medskip

{\bf Main Theorem }{\it 
Let $Q$ be a quadric bundle on a space $X$, generically non-degenerate 
of rank $r \ge 2$, 
mildly degenerating over a topological divisor $Y\subset X$.
Let $Q_{X-Y}$ be the resulting rank $r$ non-degenerate quadric bundle on
$X-Y$, and let $\T^{Q_Y}$ be the resulting rank $r-1$ non-degenerate 
canonical triple on $Y$. 
Let $\alpha\in PH^*(BGO(r))$ be a universal quadric invariant, and let 
$\alpha(Q_{X-Y})\in H^*(X-Y)$ be its value on $Q_{X-Y}$.
Let $\delta : PH^*(BGO(r))\to H^{*-1}(BGO(r-1))$ 
be the linear map defined above,
and let $(\delta(\alpha))(\T^{Q_Y})$ be the value of the resulting
$GO(r-1)$-characteristic class $\delta(\alpha)$ on 
$\T^{Q_Y}$. Let $\ovnu_Y(\det(Q)) \in H^0(Y,\Z)$ be the  
topological vanishing parity along $Y$ of the discriminant 
$\det(Q) \in \Gamma(X,L^r\otimes \det(E)^{-2})$. 
Then under the Gysin boundary map $\pa : H^*(X-Y) \to H^{*-1}(Y)$, 
we have the equality
$$\pa(\alpha(Q_{X-Y})) = \ovnu_Y(\det(Q)) \cdot (\delta(\alpha))(\T^{Q_Y})$$
}

\medskip

{\bf Sketch of the proof } If $X$ is an algebraic variety with an effective 
Cartier divisor $Y$, and $E$ is a vector bundle on $X$ with a given line
subbundle $E'$ of the restriction $E|_Y$, then recall that the 
Hecke transform $F$ of $(E,E')$ along $Y$ 
is the vector bundle (locally free sheaf) on $X$ 
consisting of germs of all sections of $E$ which when restricted to
$Y$ lie in $E'$. Doing what may be described as a topological analog of this
procedure, we first reduce to the case where the 
absolute topological vanishing multiplicity 
of $\det(Q)$ along each component of $Y$ is either $0$ or $1$, depending on
whether the original absolute topological
vanishing multiplicity was even or odd, respectively. 

Without loss of generality, one can assume that $Y$ is connected. 
When the absolute topological vanishing multiplicity of the 
discriminant is zero along 
the divisor $Y$, we show that topologically $Q$ admits a non-degenerate 
prolongation from $X-Y$ to all of $X$. In particular, all quadric invariants
map to zero under the Gysin boundary map. This is necessarily the case when
the real vector bundle $N$ on $Y$ (the `topological normal bundle' of
$Y$ in $X$) is not orientable, as then the absolute topological 
vanishing multiplicity of the discriminant is always zero.

So now remains the case where $|\nu|_Y(\det(Q))=1$. 
(In the algebraic or complex analytic category, if the base $X$ and the
divisor $Y$ are non-singular varieties, then the condition
$\nu_Y(\det(b))=1$ on a mildly degenerating algebraic or holomorphic 
triple $(E,L,b)$ is equivalent to demanding that the closed subvariety
$V \subset P(E)$ defined by $b$ is non-singular.) 
To study the case $|\nu|_Y(\det(Q))=1$, we construct a 
particular such quadric bundle $\QQ$ which is quasi-universal in 
the sense that for any $Q$ on a base $X$, degenerating over $Y$,
there is a tubular neighbourhood of $Y$ in $X$ 
over which the original quadric bundle
may be replaced by
a pullback of the quasi-universal bundle $\QQ$, for the 
purpose of the main theorem. 
By this device, we reduce the problem to understanding 
the Gysin boundary map
$\pa$ in the case of the quasi-universal quadric bundle $\QQ$.

The quasi-universal bundle $\QQ$ is constructed as follows, 
for all $r\ge 2$. 
Let $(\EE,\LL,\bb)$ be the universal triple over $BGO(r-1)$.
The total space $\LL$ of the line bundle $\pi : \LL \to BGO(r-1)$
serves as the base space for $\QQ$, and the zero section
$BGO(r-1)\subset \LL$ is the degeneration locus of $\QQ$. 
Let $\tau \in \Gamma(\LL,\pi^*(\LL))$ be the tautological section
of the pullback of $\LL$ under $\pi$, which vanishes with 
multiplicity $1$ along the zero section $BGO(r-1)\subset \LL$. 
We may regard $\tau$ as defining a bilinear form 
$\tau : \O_{\LL}\otimes \O_{\LL} \to \pi^*(\LL)$ on the trivial line 
bundle $\O_{\LL}$ on $\LL$. The direct sum triple 
$(\O_{\LL},\pi^*(\LL),\tau)\oplus \pi^*(\EE,\LL,\bb)$ 
on $\LL$ defines the quadric bundle
$$\QQ=[\O_{\LL}\oplus\pi^*(\EE),\,\pi^*(\LL),\,\tau\oplus\pi^*(\bb)]$$
which is the desired quasi-universal mildly 
degenerating quadric bundle with generic rank $r$.
By the construction of the quasi-universal triple, the theorem follows.

The above theorem tells us how to compute explicitly $\pa(\alpha(Q))$,
once we know how to compute the Gysin boundary map 
$H^*(BO(r-1)) \to H^{*-1}(BGO(r-1))$. We write this Gysin boundary map 
explicitly in terms of generators of the cohomologies. 
The answer naturally falls into two cases.

{\bf (1) Odd rank degenerating to even rank } 
When $Q$ is non-degenerate of rank $2n+1$ on $X-Y$, degenerating to
rank $2n$ on $Y$, the required Gysin boundary map
$d : H^*(BO(2n)) \to H^{*-1}(BGO(2n))$ is explicitly given in Section 4.

{\bf (2) Even rank degenerating to odd rank } 
When $Q$ is non-degenerate of rank $2n+2$ on $X-Y$, degenerating to
rank $2n+1$ on $Y$, the required Gysin boundary map
$d : H^*(BO(2n+1)) \to H^{*-1}(BGO(2n+1))$ is already known, 
as recalled in Section 3.

The ``rank $3$ degenerating to rank $2$'' case of the above  
theorem, where a conic bundle $Q$ degenerates into a bundle whose 
fiber is a pair of distinct lines, was considered earlier in [N]. 
In this case, $PH^*(BGO(3)) = \FF_2[\wh{w}_2,\wh{w}_3]$, while
$H^*(BGO(2)) = \FF_2[\lambda, a_1, b_4]/(\lambda a_1)$ 
(as proved in [H-N]). 
The invariant $\wh{w}_2(Q)\in H^2(X-Y)$ is the Brauer class 
of the $\PP^1$-bundle $Q_{X-Y}$ on $X-Y$. The restriction 
$Q_Y$ defines a $2$-sheeted cover of $Y$, and the quadric invariant 
$a_1 \in H^1(Y)$ is the class of this cover.

It was proved in [N] that $\pa(\wh{w}_2) = \nu_Y(\det(Q)) a_1 $,  
but nothing was proved there about the behavior of 
the general invariant $f(\wh{w}_2,\wh{w}_3)$ under $\pa$.
As a consequence of our Main Theorem, we now know 
$\pa(f(\wh{w}_2,\wh{w}_3))$ explicitly. 
For example,  
$\wh{w}_3 \mapsto \nu_Y(\det(Q)) a_1^2$, and 
$\wh{w}_2^3 \mapsto \nu_Y(\det(Q))(a_1^5 + a_1b_4)$.
This last example is noteworthy, as $a_1^5 + a_1b_4$ is a characteristic 
class for $GO(2)$ which is not a quadric invariant.

This paper is arranged as follows.
The Sections 2 to 4 deal with non-degenerate quadric bundles.
In Section 2, the cohomological invariants of non-degenerate 
quadric bundles in 
odd ranks (which is the easy case) are
described, together with the Gysin boundary map 
$d : H^*(BO(2n+1)) \to H^{*-1}(BGO(2n+1))$.  
In Section 3, which treats even ranks $r=2n$, 
we first describe (following [H-N])
the characteristic classes for
$GO(2n)$ and the Gysin boundary map $d : H^*(BO(2n)) \to H^{*-1}(BGO(2n))$. 
Then we determine the the ring homomorphism 
$H^*(BGL(2n))\to H^*(BGO(2n))$ induced by the inclusion 
$GO(2n)\hra GL(2n)$, which should be of independent interest.
This in particular allows us to write the action of 
$\T \mapsto \T \otimes K$ 
on the $GO(2n)$-characteristic classes, which gives a 
computational procedure
to decide whether a given characteristic class is quadric invariant.
The ring of quadric invariants in ranks $4$ and $4m+2$ is described in
the Section 4, making use of the corresponding results of Toda [T] for
the orthogonal groups.  
In Section 5, we recall the properties of the Gysin boundary map that we 
need, and prove a basic lemma about the Gysin boundary map. 
In Section 6, we consider the topological behavior of
sections of complex line bundles on a space $X$ which vanish on a 
topological divisor $Y$.
In Section 7, we establish some
basic properties of mildly degenerating 
triples, and complete the proof our main theorem for 
Gysin boundary behavior of quadric invariants.
In Section 8, we show how the main theorem leads to an algorithm for 
calculations, and as illustrations we give 
explicitly the behavior of some of the quadric invariants
for the degenerations from ranks $3$ to $2$, $4$ to $3$, 
$5$ to $4$, and $6$ to $5$.
The paper ends with an Appendix (Section 9) which recalls the 
required results of Toda. 

\bigskip

{\bf Acknowledgments } The authors thank the referee for pointing out
that as the proof of the Main Theorem does not make essential use of 
any hypothesis about coefficients, the statement should be formulated 
in this generality (instead of just for coefficients $\FF_2$ as in the 
earlier version). The authors respectively thank 
the Abdus Salam ICTP Trieste, and the University of Essen, for support 
while part of this work was done.

%

\section{Quadric invariants in odd ranks}

\centerline{\bf Some generalities for all ranks}

All topological spaces will be assumed to be Hausdorff, and 
of the homotopy type of a CW complex. In this section we assume that  
all cohomologies will be singular cohomology with
coefficients in the field $\FF_2=\Z/(2)$, unless otherwise indicated.

By a Gram-Schmidt argument, isomorphism classes of 
principal $O(r)$-bundles on a base 
$X$ are the same as isomorphism classes of pairs $(E,q)$ where 
$E$ is a vector bundle on $X$
of rank $r$, and $q: E\otimes E \to \O_X$ is an everywhere non-degenerate 
symmetric bilinear form. More generally, 
isomorphism classes of principal $GO(r)$-bundles on a base 
$X$ are the same as isomorphism classes of triples $(E,L,b)$ 
where $E$ is a vector bundle on $X$
of rank $r$, $L$ is a line bundle on $X$, 
and $b: E\otimes E \to L$ is an everywhere non-degenerate 
symmetric bilinear form. It follows that isomorphisms 
$L\to \O_X$ are the same as reductions of structure group 
from $GO(r)$ to $O(r)$. In particular, if $L$ denotes the total space of
the line bundle $\pi : L\to X$, then the pullback of 
$(E,L,b)$ to $L_o = L - X$ (complement of the zero section of $L$)
has a canonical reduction of structure group to  $O(r)$. 

Let the homomorphism $\sigma : GO(r)\to \C\,^*$ be
defined by ${^t}gg =\sigma(g)I$. This fits in the short exact sequence
$$1 \to O(r)\to GO(r) \stackrel{\sigma}{\to} \C\,^* \to 1$$
Hence the classifying space $BO(r)$ is the principal $\C\,^*$-bundle on 
$BGO(r)$ associated to $\sigma$.
Let $(\EE,\LL,\bb)$ be the universal triple on the classifying space
$BGO(r)$. Then $\LL$ is associated to $\sigma$, which (again) shows that 
the complement $\LL_o$  of the zero section of $\LL$ may be identified with
$BO(r)$.

Let $\Z/(2) = \{ \pm 1\}$, and let $\mu : \Z/(2) \times O(r) \to O(r)$
be the multiplication map, sending $(a,g) \mapsto ag$. Similarly,
let $\mu : \C\,^* \times GO(r) \to GO(r)$ be 
the multiplication map, sending $(a,g) \mapsto ag$. 
As these maps are group homomorphisms, which
commute with the inclusion homomorphisms $\Z/(2)\hra \C\,^*$ and
$O(r)\hra GO(r)$, we have the following commutative diagram of
classifying spaces, 
where the vertical maps are induced by the inclusions.
$$\begin{array}{ccc}
B(\Z/(2)) \times BO(r) & \stackrel{B(\mu)}{\to} &  BO(r)  \\
\dna                   &                        &  \dna   \\
B\C\,^* \times BGO(r)  & \stackrel{B(\mu)}{\to} &  BGO(r)
\end{array}$$

Recall that $B(\Z/(2)) = \R\PP^{\infty}$, and its
cohomology ring is the polynomial ring $\FF_2[w]$  
where $w\in H^1(B(\Z/(2)))$. Also, 
$B\C\,^* = \C\PP^{\infty}$, and its 
cohomology ring is the polynomial ring $\FF_2[t]$ 
where $t\in H^2(B\C\,^*)$ is the 
image under $H^2(B\C\,^*;\Z)\to H^2(B\C\,^*)$ of the first  
Chern class of the line bundle 
$\O(1)$ on $\C\PP^{\infty}$.
The inclusion $\Z/(2) \hra \C\,^*$ induces the ring homomorphism
$H^*(B\C\,^*)\to H^*(B(\Z/(2)))$ under which 
$c$ maps to $w^2$.
Hence we get the following.

\rem\label{quad inv to ortho quad inv} 
The following  diagram of ring homomorphisms is commutative
where the first vertical map $\pi^*$ is induced by the projection
$\pi : BO(r) \to BGO(r)$, and the second vertical map $\theta$ 
is induced by $\pi^*$ together with $w \mapsto t^2$. 
$$\begin{array}{ccc}
H^*(BGO(r)) & \stackrel{B(\mu)^*}{\lra} & H^*(BGO(r))[t]  \\
{\scriptstyle \pi^*}\dna ~~ & & ~~ \dna {\scriptstyle \theta}\\
H^*(BO(r)) & \stackrel{B(\mu)^*}{\lra} & H^*(BO(r))[w]
\end{array}$$

In the rest of this section, we record various facts about
non-degenerate triples and quadric bundles in odd ranks, for later use.

\bigskip

\centerline{\bf The cohomology ring $H^*(BGO(2n+1))$}

The map $\C\,^*\times SO(2n+1)\to GO(2n+1)$ defined by
$(a,h)\mapsto ah$ is an isomorphism. 
As $H^*(B\C\,^*) = \FF_2[c]$ 
where $c$ has degree $2$, while 
$H^*(BSO(2n+1))= \FF_2[\wh{w}_2,\ldots,\wh{w}_{2n+1}]$
where $\wh{w}_i \in H^i(BSO(2n+1))$ 
are the universal special Stiefel-Whitney classes for $SO(2n+1)$,
we get an identification
$$H^*(BGO(2n+1)) = \FF_2[c,\wh{w}_2,\ldots,\wh{w}_{2n+1}]$$
where we write $c$ for $t\otimes 1$, and write $\wh{w}_i$ for
$1\otimes \wh{w}_i$.

\medskip

\centerline{\bf Quadric invariants in odd rank} 

\begin{proposition}\label{projoddinv} 
The ring of quadric invariants in any odd rank $2n+1$ is the subring 
$H^*(BSO(2n+1))= \FF_2[\wh{w}_2,\,\ldots,\,\wh{w}_{2n+1}]$
of $H^*(BGO(2n+1))$.
\end{proposition}


\medskip

\centerline{\bf Gysin sequence for $BO(1)\to B\C\,^*$}

The character $2\chi : \C\,^*\to \C\,^* : z\mapsto z^2$ has kernel 
$O(1) = \{\pm 1\}\subset \C\,^*$, hence the classifying space
$BO(1)$ can be taken to be the complement of the zero section
of the line bundle $\O(2)$ on $B\C\,^* = \C\PP^{\infty}$ 
defined by the character $2\chi$.
The Euler class of $\O(2)$
is $0\in H^2(B\C\,^*)$ as $2=0$ in $\FF_2$. 
Hence the long exact Gysin sequence for the principal $\C\,^*$-bundle
$\pi : BO(1)\to B\C\,^*$ splits to give short exact sequences 
$$0\to H^i(B\C\,^*) \stackrel{\pi^*}{\to} H^i(BO(1))
\stackrel{d}{\to} H^{i-1}(B\C\,^*)\to 0$$
For all $i\ge 0$, 
$H^{2i}(B\C\,^*) = \{0,c^i\}$, $H^{2i+1}(B\C\,^*) = 0$,
and $H^i(BO(1)) = \{0, w^i\}$, 
where we take $c^0=1$ and $w^0=1$.
Hence the above short exact Gysin sequences
implies that
$$\pi^*(c^i) = w^{2i},~~d(w^{2i}) = 0, \mbox{ and }
d(w^{2i+1}) = c^i, \mbox{ for all } i\ge 0.$$

\medskip

\centerline{\bf Gysin sequence for $BO(2n+1) \to BGO(2n+1)$}

The classifying space
$BO(2n+1)$  has cohomology ring
$\FF_2[w_1,\,\ldots,\,w_{2n+1}]$,
the polynomial ring in the Stiefel-Whitney classes $w_i$.
With $O(1) = \{\pm 1\} \subset \C\,^*$, we
have an isomorphism $O(1) \times SO(2n+1)
\to O(2n+1) : (a, h)\mapsto ah$, which gives an isomorphism
$BO(1) \times BSO(2n+1) \to BO(2n+1)$. In our previous notation,
$H^*(BO(1)) = \FF_2[w]$ and $H^*(BSO(2n+1)) = 
\FF_2[\wh{w}_2,\,\ldots,\,\wh{w}_{2n+1}]$, so taking tensor product,
$$H^*(BO(2n+1)) = \FF_2[w,\,\wh{w}_2,\,\ldots,\,\wh{w}_{2n+1}]$$
Note that if $s_1,\ldots,s_m$
are the elementary symmetric functions in the variables
$x_1,\ldots,x_m$, and $u$ is another variable, then we have
$$s_r(u+x_1,\ldots,u+x_m) = \sum_{0 \le i \le r} 
{m-i \choose r-i}\, u^{r-i}\, s_i(x_1,\ldots,x_m)$$
Hence by the splitting principle we get the following.

\begin{lemma}\label{w in terms of wbar} 
In $H^*(BO(2n+1)) = \FF_2[w_1,\,\ldots,\,w_{2n+1}]
=\FF_2[w,\,\wh{w}_2,\,\ldots,\,\wh{w}_{2n+1}]$, we have the identities
\begin{eqnarray*} 
w_1 & = & w, \\ 
w_r & = & {2n+1 \choose r}\,w^r + \sum_{2\le i \le r}   
{2n+1-i \choose r-i}\, w^{r-i}\, \wh{w}_i ~~\mbox{for } 2\le r \le 2n+1,\\
\wh{w}_r & = & {2n \choose r}\,w_1^r + \sum_{2\le i \le r} 
{2n+1-i \choose   r-i}\, w_1^{r-i}\,w_i  ~~\mbox{for } 2\le r \le 2n+1. 
\end{eqnarray*} 
\end{lemma}

\begin{lemma}\label{Gysin for BO(2n+1) to BGO(2n+1)} 
For the $\C\,^*$-fibration $\pi : BO(2n+1)\to BGO(2n+1)$, 
we have the following.

(i) The Euler class is zero.

(ii) The ring homomorphism 
$\pi^* : H^*(BGO(2n+1))\to H^*(BO(2n+1))$ is given in terms of 
generators by
\begin{eqnarray*} \pi^*(c) & = & w_1^2, \\ 
\pi^*(\wh{w}_r) & = & {2n\choose r}\,w_1^r 
+ \sum_{2\le i \le r} {2n+1-i \choose  r-i}\, 
w_1^{r-i}\, w_i  ~~\mbox{for } 2\le r \le 2n+1. 
\end{eqnarray*}  

(iii) The Gysin boundary map $d : H^*(BO(2n+1)) \to  
H^{*-1}(BGO(2n+1))$ has the following expression. Using the identities
given by Lemma \ref{w in terms of wbar}, any element of 
$H^*(BO(2n+1))=\FF_2[w_1,\ldots,w_{2n+1}]$ 
can be uniquely expressed as a polynomial
$\sum w^i f_i(\wh{w}_2,\ldots,\wh{w}_{2n+1})$,
where the $f_i$ are polynomials in $2n$ variables.
Then the additive group homomorphism $d$ is given by 
$$d( \sum_{i\ge 0} w^i f_i(\wh{w}_2,\ldots,\wh{w}_{2n+1})) 
= \sum_{j\ge 0} c^j f_{2j+1}(\wh{w}_2,\ldots,\wh{w}_{2n+1})$$
\end{lemma} 

\proof The isomorphisms $O(1) \times SO(2n+1) 
\to O(2n+1) : (a, h)\mapsto ah$ and
$\C\,^* \times SO(2n+1)
\to GO(2n+1) : (a, h)\mapsto ah$ fit in the commutative square
$$\begin{array}{ccc}
O(1) \times SO(2n+1) & \to & O(2n+1) \\
\dna                &     & \dna    \\
\C\,^* \times SO(2n+1) & \to & GO(2n+1)
\end{array}$$
Hence by the previous calculation of the  
Gysin sequence for $BO(1)\to B\C\,^*$, the 
Gysin sequence for $BO(2n+1) \to BGO(2n+1)$ 
breaks into short exact sequences
$$0\to H^i(BGO(2n+1)) \stackrel{\pi^*}{\to} 
H^i(BO(2n+1)) \stackrel{d}{\to} H^{i-1}(BGO(2n+1))\to 0$$
with $\pi^*(c) = w^2$, 
$\pi^*(\wh{w}_r) = \wh{w}_r$ for all $r\ge 2$,
$d(w^{2j}) = 0$ for all $j\ge 0$, 
$d(w^{2j+1}) = c^j$ for all $r\ge 0$, and
$d(f(\wh{w}_2,\ldots,\wh{w}_{2n+1})) = 0$ for any
polynomial $f$ in the variables $\wh{w}_2,\ldots,\wh{w}_{2n+1}$.

As $f(\wh{w}_2,\ldots,\wh{w}_{2n+1})$ is pulled back from
$BGO(2n+1)$, it follows from
Lemma \ref{cupgysin} that 
$d(w^if(\wh{w}_2,\ldots,\wh{w}_{2n+1})) = 
d(w^i)f(\wh{w}_2,\ldots,\wh{w}_{2n+1})$, 
which completes the proof.
\hfill$\square$

\example 
In particular $d(w_1) = 1$, and for all $r\ge 1$ 
we have $d(w_{2r})=0$ and $d(w_{2r+1}) = 
{2n \choose 2r}\, c^r + 
\sum_{1\le j \le r}{2n-2j \choose 2r-2j}\, c^{r-j}\,\wh{w}_{2j}$.

%
%

\rem 
The ring $H^*(BGL(2n+1);\Z)$ equals $\Z[c_1,\ldots,c_{2n+1}]$
where the $c_i$ are the Chern classes, therefore modulo $2$ we
have $H^*(BGL(2n+1))=\FF_2[\ov{c}_1,\ldots,\ov{c}_{2n+1}]$
with $\ov{c}_i \in H^{2i}(BGL(2n+1))$.
Under the map $H^*(BGL(2n+1)) \to H^*(BGO(2n+1))$, it can be seen that
$$\ov{c}_r \mapsto {2n+1\choose  r}c^r + 
\sum_{2\le i\le r}{m-i\choose  r-i}c^{r-i}\wh{w}_i^2$$

%

\section{Characteristic classes for $GO(2n)$}

\centerline{\bf Generators $\lambda$,
$a_{2i-1}$, $b_{4j}$, $d_T$ of $H^*(BGO(2n))$, and relations.}

The cohomology ring $H^*(BGO(2n))$ with coefficients $\FF_2$ 
has been explicitly determined in terms of generators and relations
in [H-N], which we recall. Let $(\EE,\LL,\bb)$ denote the
universal triple on $BGO(2n)$. 
The characteristic class $\lambda \in H^2(BGO(2n))$ is by definition
the Euler class of $\LL$. 

The character $\sigma :GO(2n)\to \C\,^*$
defined by ${^t}gg = \sigma(g)I$ has kernel $O(2n)$, 
and $\LL$ is associated to the character $\sigma$, hence
$BO(2n)$ can be taken to be the complement $\LL_o$ of the zero section
of the line bundle $\LL$. For the principal $\C\,^*$-bundle
$\pi : BO(2n)\to BGO(2n)$, consider the long exact Gysin sequence
{\small
$$\cdots \stackrel{\lambda}{\to} H^i(BGO(2n)) 
\stackrel{\pi^*}{\to} H^i(BO(2n)) 
\stackrel{d}{\to} H^{i-1}(BGO(2n)) 
\stackrel{\lambda}{\to} H^{i+1}(BGO(2n)) 
\stackrel{\pi^*}{\to} \cdots$$
}
For each $1\le j\le n$, we put 
$$a_{2j-1} = dw_{2j}\in H^{2j-1}(BGO(2n)).$$
More generally, for any subset 
$T=\{ i_1,\ldots,i_r\} \subset \{ 1,\ldots,n\}$ 
of cardinality $r\ge 2$, let 
$v_T = w_{2i_1}\cdots w_{2i_r} \in H^{2\deg(T)}(BO(2n))$
where $\deg(T) = i_1 + \cdots + i_r$. We put
$$d_T = d(v_T) \in H^{2\deg(T) - 1}(BGO(2n)).$$
We recall from [H-N] that the composite 
$\pi^*d : H^*(BO(2n)) \to H^{*-1}(BO(2n))$ is the derivation
$s = \sum w_{2i-1}{\pa \over \pa w_{2i}}$. Hence we have
$$\pi^*(a_{2j-1}) = w_{2j-1} ~~\mbox{and }
\pi^*(d_T) = \sum_{i\in T} w_{2i-1}v_{T-\{ i\}}
\mbox{ in }  H^*(BO(2n)).$$
Finally, let $c_i(\EE)\in H^{2i}(BGL(2n);\Z)$ denote the $i$ th Chern
class of the vector bundle $\EE$, and let  
$\ov{c}_i(\EE)\in H^{2i}(BGL(2n))$ be its image under 
the coefficient map $\Z\to \FF_2$. 
We put 
$$b_{4j} = \ov{c}_{2j}(\EE) \in H^{4j}(BGO(2n))
~~\mbox{for each } 1\le j\le n.$$
These satisfy
$$\pi^*(b_{4j}) = w_{2j}^2~~\mbox{for each } 1\le j\le n.$$

As shown in Theorem 3.3.5 of [H-N], these elements 
$\lambda, a_{2i-1}$, $b_{4j}$, and $d_T$ generate
the ring $H^*(BGO(2n))$, with relations given as follows.

\begin{theorem}\label{cohomology of BGO(2n)} 
For any $n\ge 1$, the cohomology ring of $BGO(2n)$ with 
coefficients $\FF_2$ is the quotient 
$$H^*(BGO(2n))= 
{\FF_2[\lambda, (a_{2i-1})_i, (b_{4i})_i, (d_T)_T ]\over  I}$$
where $I$ is the ideal generated by
the elements $\lambda a_{2i-1}$ for $1\le i\le n$, the 
elements $\lambda d_T$ 
where $T \subset \{ 1,\ldots,n\}$ is a subset of cardinality 
$|T|\ge 2$, 
the elements 
$\sum_{i\in T} a_{2i-1} \,d_{T-\{ i\}}$ 
where $T \subset \{ 1,\ldots,n\}$ with $|T| \ge 3$, 
the elements
$(d_{\{ i,j \} })^2 + a_{2i-1}^2b_{4j} + a_{2j-1}^2b_{4i}$ 
where $\{ i,j\}\subset \{ 1,\ldots,n\}$ is 
a subset of cardinality $2$, and 
the elements 
$d_T\,d_U + \sum_{p\in T} \prod_{q \in T \cap U -\{ p\} }  
\, a_{2p-1} \, b_{4q} \, d_{(T - \{ p\}) \Delta U}$
where $T\ne U$ in case $|T|=|U|=2$,
and where $\Delta$ denotes the symmetric difference of sets.
\end{theorem}

\centerline{\bf The map $\pi^*: H^*(BGO(2n)) \to H^*(BO(2n))$}

The pullback map $\pi^*: H^*(BGO(2n)) \to H^*(BO(2n))$
satisfies the following equations, as shown in [H-N]. 
\begin{eqnarray*} 
\pi^*(\lambda) & = & 0,\\
\pi^*(a_{2i-1}) & = & w_{2i-1} ~~\mbox{for all } 1\le i\le n,\\
\pi^*(b_{4i}) & = & w_{2i}^2 ~~\mbox{for all } 1\le i\le n,\\
\pi^*(d_T) & = & \sum_{i\in T}w_{2i-1}v_{T-\{ i\}}  ~~\mbox{for all } 
             T \subset \{ 1,\ldots,n\} \mbox{ with } |T|\ge 2, 
\end{eqnarray*} 
where $v_J = \prod_{j\in J}w_{2j}$
for any subset $J \subset \{ 1,\ldots,n\}$.

Moreover, the map 
$\pi^* : H^{2i+1}(BGO(2n)) \to H^{2i+1}(BO(2n))$
is injective for each $i\ge 0$, by Remark 3.3.4 of [H-N].
Hence the Gysin sequence breaks into exact sequences
$$H^{2i}(BGO(2n))\stackrel{\pi^*}{\to}H^{2i}(BO(2n))
\stackrel{d}{\to}H^{2i-1}(BGO(2n))\to 0$$
and
$$0\to H^{2i+1}(BGO(2n))\stackrel{\pi^*}{\to}H^{2i+1}(BO(2n))
\stackrel{d}{\to}H^{2i}(BGO(2n))$$

\centerline{\bf Gysin boundary map $d : H^*(BO(2n))\to H^{*-1}(BGO(2n))$} 

The Gysin boundary map 
$d : H^*(BO(2n))\to H^{*-1}(BGO(2n))$ satisfies the following 
equations, as shown in [H-N].  
\begin{eqnarray*}
d(w_{2i-1}) & = & 0 ~~\mbox{for all } 1\le i\le n,~\mbox{as }
                   w_{2i-1} = \pi^*(a_{2i-1}),\\
d(w_{2i}^2) & = & 0 ~~\mbox{for all } 1\le i\le n,~\mbox{as }
                    w_{2i}^2 = \pi^*(b_{4i}),\\
d(w_{2i}) & = & a_{2i-1} ~~\mbox{for all } 1\le i\le n,\\
d(w_{2i_1}\cdots w_{2i_r}) & = & d_T ~~\mbox{for all } 
             T = \{ i_1 <\ldots < i_r \}\subset \{ 1,\ldots,n\}. 
\end{eqnarray*}
Moreover, by Lemma \ref{cupgysin} proved later, it follows
that for any $x\in H^*(BGO(2n))$ and $y\in H^*(BO(2n))$ we have
$d(\pi^*(x)y) = \pi^*(x)dy$. Hence $d : H^*(BO(2n))\to H^{*-1}(BO(2n))$ 
is determined for any monomial in the $w_i$'s as follows.
Any such monomial can be uniquely written  
as a product 
$f(w_{2i-1})\cdot g(w_{2i}^2)\cdot w_{2i_1} \cdots w_{2i_r}$ where 
$f$ and $g$ are monomials in $n$ variables, and 
$1\le i_1 < i_2 < \ldots < i_r \le n$.
Then we have
$$d(f(w_{2i-1})\cdot g(w_{2i}^2)\cdot w_{2i_1} \cdots w_{2i_r}) =
\left\{ 
\begin{array}{ll}
0 & \mbox{if } r= 0 \\
f(a_{2i-1})\cdot g(b_{4i})
\cdot\, a_{2i_1-1} & \mbox{if } r=1 \\
f(a_{2i-1})\cdot g(b_{4i})
\cdot d_T & \mbox{if } r \ge 2, ~\mbox{where } \\
& T = \{ i_1 <\ldots < i_r \}.
\end{array}\right.
$$
This completes our description of the Gysin boundary map
for $BO(2n) \to BGO(2n)$.

\medskip

\centerline{\bf The map $H^*(BGL(2n))\to H^*(BGO(2n))$}

The determination of the ring homomorphism 
$H^*(BGL(2n))\to H^*(BGO(2n))$, which is carried out in the rest of
this section, is needed in order to determine the 
quadric invariants. 

The cohomology ring $H^*(BGL(2n))$ is the polynomial ring
$\FF_2[\ov{c}_1, \ldots, \ov{c}_{2n}]$ where the $2n$ variables
$\ov{c}_r$ are the Chern classes mod $2$.
By the definition of the elements 
$b_{4r} \in H^{4r}(BGO(2n))$ in terms of the  
universal triple $(\EE,\LL,\bb)$,
we have 
$\ov{c}_{2r}(\EE) = b_{4r}$ for $1\le i\le n$. 
Hence under the ring homomorphism
$H^*(BGL(2n))\to H^*(BGO(2n))$ induced by the inclusion 
$GO(2n)\hra GL(2n)$, we have 
$\ov{c}_{2r}\mapsto b_{4r}$. The following proposition gives
the images of the odd Chern classes $\ov{c}_{2r-1}$ in $H^*(BGO(2n))$,
which completes the determination of the ring homomorphism 
$H^*(BGL(2n))\to H^*(BGO(2n))$.

\begin{proposition}\label{odd Chern classes} 
Consider the lower triangular $n\times n$ matrix $A$ over 
$\FF_2[\lambda^2]$, 
with all diagonal entries $1$, and below-diagonal entries given by
$$A_{r,k} = {n-k \choose 2r - 2k} \lambda^{2r-2k} \mbox{ for } r>k.$$
Let $B = A^{-1}$ be its matrix inverse over the ring $\FF_2[\lambda^2]$.
For $1\le r\le n$, let the polynomial
$f_{n,r}(\lambda, b_4,\ldots, b_{4r-4})$ be defined by
{\small 
$$f_{n,r} = {n \choose 2r-1} \lambda^{2r-2} 
+\sum_{1\le k\le r-1}{n-k \choose 2r-1-2k} \lambda^{2r-2-2k}
\left(\sum_{1\le j\le k}B_{k,j}(b_{4j} - 
{n \choose 2j} \lambda^{2j})\right)$$
} 

Then in the cohomology ring $H^*(BGO(2n))$, we have 
$$\ov{c}_{2r-1}(\EE) =  
a_{2r-1}^2  + \lambda \cdot f_{n,r}(\lambda, b_4,\ldots, b_{4r-4})$$
\end{proposition}

\proof We divide the proof into two steps, {\bf (a)} and  {\bf (b)}.

{\bf (a) } For each $1\le r\le n$, there exists a unique polynomial 
$g_r(\lambda, b_4,\ldots, b_{4r-4})$ such that
$\ov{c}_{2r-1}(\EE) =  
a_{2r-1}^2  + \lambda \cdot g_r(\lambda, b_4,\ldots, b_{4r-4})$.

{\bf (b) } The polynomials 
$g_r(\lambda, b_4,\ldots, b_{4r-4})$ are the polynomials 
$f_{n,r}(\lambda, b_4,\ldots, b_{4r-4})$ defined 
in the statement of the proposition.

{\bf Proof of (a) } The composite homomorphism
$O(2n)\hra GO(2n)\hra GL(2n)$
induces $H^*(BGL(2n))\to H^*(BGO(2n))\to H^*(BO(2n))$
under which $\ov{c}_i \mapsto w_i^2$.
As $a_{2i-1}^2\mapsto w_{2i-1}^2$ under $H^*(BGO(2n))\to H^*(BO(2n))$,
we must have 
$$\ov{c}_{2i-1}(\EE) = a_{2i-1}^2 + z_{4i-2}$$ 
where $z_{4i-2}$ lies in the kernel of 
$\pi^*: H^{4i-2}(BGO(2n))\to H^{4i-2}(BO(2n))$.
By exactness of Gysin, we have $z_{4i-2} = \lambda y_{4i-4}$
for some $y_{4i-4} \in H^{4i-4}(BGO(2n))$. 
Now, by the structure of the ring
$H^*(BGO(2n))$, we know that $\lambda$ annihilates 
the $a_{2j-1}$'s and the $d_T$'s. Hence we can replace 
$y_{4i-4}$ by a polynomial $g_i(\lambda, b_{4j})$ 
in $\lambda $ and the $b_{4j}$'s. As the variables $\lambda$ and 
the $b_{4j}$ are algebraically independent, $g_i$ is unique. 
By degree considerations, the highest $b_{4j}$ that occurs 
in $g_i$ can be at most $b_{4i-4}$. This completes the proof of (a).

{\bf Proof of (b) } We first determine 
the polynomial $g_i(\lambda, b_4,\ldots, b_{4i-4})$ in the following
Example \ref{lambda and b} of a triple $\T$. 

\example\label{lambda and b} 
Let $\O(1)$ be the universal line bundle on $B\C\,^* = \PP^{\infty}_{\C}$.
Let $X = B\C\,^*\times \ldots \times B\C\,^*$ 
be the product of $n+1$-copies,
and let $p_i : X \to B\C\,^*$ be the projections for $0\le i \le n$.
Let $L = p_0^*(\O(1))$, and for $1\le i\le n$ let
$K_i = p_i^*(\O(1))$. On $X$, we get the non-degenerate triples
$\T_i = ((L\otimes K_i)\oplus K_i^{-1}, L, b_i)$ where $b_i$ is induced by
the canonical isomorphism 
$(L\otimes K_i) \otimes K_i^{-1} \stackrel{\sim}{\to} L$. 
Let $\T = (E,L,b)$ be the direct sum triple $\oplus \T_i$
(see Definition \ref{triple sum}). 
We now write down the characteristic classes of $\T$. 
Let $\lambda = \ov{c}_1(L)$. Then by definition, 
$\lambda(\T)= \lambda$. 
As the odd cohomologies of $X$ are zero, the classes $a_{2i-1}(\T)$
and $d_T(\T)$ are zero for all $1\le i \le n$ and for all
$T\subset \{ 1,\ldots, n\}$ with $|T| \ge 2$. 

Let $x_i = \ov{c}_1(K_i)$ for $1\le i\le n$, and let 
$s_k(\lambda x_i+ x_i^2)$ denote the $k$ th elementary symmetric
polynomial in the variables 
$\lambda x_1+ x_1^2,\ldots, \lambda x_n+ x_n^2$.
As $E$ is the direct sum of the $(L\otimes K_i) \oplus K_i^{-1}$,
its Chern classes (modulo $2$) are given by 
\begin{eqnarray*}
\ov{c}_{2r}(E) & = & {n \choose 2r}\lambda^{2r} +
                 \sum_{1\le k\le r}{n-k \choose 2r-2k} \lambda^{2r-2k}
                 s_k(\lambda x_i+ x_i^2) \\
\ov{c}_{2r-1}(E) & = & {n \choose 2r-1} \lambda^{2r-1} +
\sum_{1\le k\le r-1}{n-k \choose 2r-1-2k} \lambda^{2r-1-2k}
                 s_k(\lambda x_i+ x_i^2) 
\end{eqnarray*}
As $b_{4r}(\T) = \ov{c}_{2r}(E)$, we have the following
equations for $1\le r\le n$.
$$b_{4r} -{n \choose 2r}\lambda^{2r}  
= \sum_{1\le k\le r}{n-k \choose 2r-2k} \lambda^{2r-2k} 
s_k(\lambda x_i + x_i^2)$$
This is a system of linear equations with coefficients in $\FF_2[\lambda^2]$, 
for the $b_{4r}-{n \choose 2r} \lambda^{2r}$ in terms of the 
$s_k(\lambda x_i + x_i^2)$. 
It is given by the $n\times n$ matrix $A$ over $\FF_2[\lambda^2]$, 
with entries
$$A_{r,k} = {n-k \choose 2r - 2k} \lambda^{2r-2k}$$
This is a lower triangular matrix, with all diagonal entries equal to $1$,
so is invertible over the ring $\FF_2[\lambda^2]$.
Let $B = A^{-1}$ be its matrix inverse. 
Hence we get
$$s_r(\lambda x_i + x_i^2)  = \sum_{1\le k\le r}B_{r,k}
(b_{4k} - {n \choose 2k} \lambda^{2k})$$
Substituting this in the equation for 
$\ov{c}_{2r-1}(E)$, we get equations 
$$ \ov{c}_{2r-1}(E) = 
\lambda \cdot f_{n,r}(\lambda, b_4,\ldots, b_{4r-4})$$
where 
$$f_{n,r} = {n \choose 2r-1}\lambda^{2r-2} +
\sum_{k=1}^{r-1}{n-k \choose 2r-1-2k}
             \lambda^{2r-2-2k}\left(
\sum_{1\le j\le k}B_{k,j}
(b_{4j} - {n \choose 2j} \lambda^{2j})\right)$$

\medskip

{\bf Proof of (b) continued } In the Example \ref{lambda and b},
we have the desired equality $g_i = f_{n,i}$.
Note that the cohomology ring $H^*(X)$ is the polynomial ring
$\FF_2[\lambda, x_1,\ldots,x_n]$, in which 
the $n+1$ elements $\lambda, b_4,\ldots, b_{4n}$ 
are algebraically independent, where 
$b_{4r} = {n \choose 2r}\lambda^{2r} +
\sum_{1\le k\le r}{n-k \choose 2r-2k} \lambda^{2r-2k} 
s_k(\lambda x_i+ x_i^2)$. Hence as
$g_i = f_{n,i}$ in this example, we get $g_i = f_{n,i}$ universally.

This completes the proof of Proposition \ref{odd Chern classes}.
\hfill$\square$

\medskip

\example 
The above proposition in particular gives the 
following identities in $H^*(BGO(2n))$.
\begin{eqnarray*}
\ov{c}_1 & = & a_1^2 + n \lambda ~~\mbox{for all } n \ge 1, ~\mbox{and}\\
\ov{c}_3 & = & a_3^2 + {n(n-1)(2n-1) \over 6} \lambda^3 + (n-1)\lambda b_4 
                    ~~\mbox{for all } n \ge 2.
\end{eqnarray*}

%

\section{Quadric invariants in even ranks}

\centerline{\bf The ring homomorphism 
$B(\mu)^* : H^*(BGO(2n))\to H^*(BGO(2n))[t]$} 

\begin{proposition}\label{action} 
The ring homomorphism 
$$B(\mu)^* : H^*(BGO(2n))\to H^*(B\C\,^*)\otimes H^*(BGO(2n)) 
= H^*(BGO(2n))[t]$$  
induced by the group homomorphism 
$\mu : \C\,^* \times GO(2n) \to GO(2n)$ which sends $(a,g)\mapsto ag$,
is given in terms of the generators of $H^*(BGO(2n))$ as follows.
\begin{eqnarray*} 
B(\mu)^*\,\lambda & = & \lambda  \\ 
B(\mu)^*\,b_4 & = & b_4 + (a_1^2 + n\lambda)t + n  t^2 \\
B(\mu)^*\,b_{4r} & = & \sum_{i=1}^{r} {2n-2i \choose 2r-2i}
( b_{4i}\, + \lambda f_{n,i} t + a_{2i-1}^2t )t^{2r-2i} 
+ {2n \choose 2r} t^{2r}\\
& & \mbox{where the elements } f_{n,i}
    \mbox{ are as in Proposition \ref{odd Chern classes}.}\\
B(\mu)^*\,a_1 & = & a_1 \\  
B(\mu)^*\,a_{2r-1} & = & \sum_{i=1}^r {2n+1-2i \choose 2r-2i} 
                            a_{2i-1}t^{r-i} \\
B(\mu)^*\,d_{\{ 1,2\}} & = & d_{\{ 1,2\}} + na_3 t + 
{n \choose 2} a_1t^2 \\
B(\mu)^*\,d_{\{ p,q\}} & = & 
\sum_{i=1}^p{2n \choose 2q}{2n -2i\choose 2p-2i}a_{2i-1} t^{p+q-i} + 
\sum_{j=1}^q{2n \choose 2p}{2n -2j\choose 2q-2j}a_{2j-1} t^{p+q-j} \\
& &~~ + \sum_{i=1}^p\sum_{j=1}^qt^{p+q-i-j}
d_{\{ i,j\}} 
~\mbox{where by convention } d_{\{ k,\ell\}} = 0 \mbox{ for } k=\ell.\\
\end{eqnarray*}
(We do not give a closed formula for $B(\mu)^*\,d_T$ for 
a general $T$, but explain in the course of the
proof how to explicitly compute $B(\mu)^*\,d_T$ for any given $T$.)
\end{proposition} 

\proof If $\T=(E,L,b)$ is a rank $2n$ triple, and $K$ any line bundle, 
then $c_1(L\otimes K^2) = c_1(L) + 2c_1(K)$ for the Chern class $c_1$. 
Hence $\ov{c}_1(L\otimes K^2) = \ov{c}_1(L)$ in $\FF_2$ coefficients. 
This gives $B(\mu)^*\,\lambda = \lambda$.  

Next we determine $B(\mu)^*\,b_{4r}$, using the expression 
$\ov{c}_{2i}(E) = b_{4i}$ and the expression for 
$\ov{c}_{2j-1}(E)$ given by the Proposition \ref{odd Chern classes}.
\begin{eqnarray*} 
& & \ov{c}_{2r}(E\otimes K) \\  
& = & \sum_{p=0}^{2r} {2n-p \choose 2r-p} \ov{c}_1(K)^{2r-p}  
                             \ov{c}_p(E)\\ 
& = &\sum_{i=0}^{r} {2n-2i \choose 2r-2i} \ov{c}_1(K)^{2r-2i}   
                             \ov{c}_{2i} (E) 
      + \sum_{j=1}^{r} {2n-2j+1 \choose 2r-2j+1} \ov{c}_1(K)^{2r-2j+1}  
                             \ov{c}_{2j -1} (E)\\ 
& = &\sum_{i=0}^{r} {2n-2i \choose 2r-2i} \ov{c}_1(K)^{2r-2i} b_{4i} + 
\sum_{j=1}^{r} {2n-2j+1 \choose 2r-2j+1} \ov{c}_1(K)^{2r-2j+1}  
      (a_{2j-1}^2 + \lambda f_{n,j}) 
\end{eqnarray*}
From this, the expression for $B(\mu)^*\,b_{4r}$ follows,
using standard binomial identities modulo $2$.

Note that for any $j\ge 0$, 
any element of $H^{2j-1}(BGO(2n))$ is a sum of terms,
each of which has an $a_{2i-1}$ or a $d_T$ as a factor, and 
so is killed by $\lambda$, so the map  
$\lambda : H^{2j-1}(BGO(2n))\to H^{2j+1}(BGO(2n))$ is identically zero,
and hence $\pi^* : H^{2j+1}(BGO(2n)) \to H^{2j+1}(BO(2n))$ is injective. 
As $t$ has even degree, it follows that the
graded ring homomorphism $\theta : H^*(BGO(2n))[t] 
\to H^*(BO(2n))[w]$ of Remark \ref{quad inv to ortho quad inv} 
is injective on each graded piece of $H^*(BGO(2n))[t]$ of odd degree. 
It therefore follows from the commutative diagram in  
Remark \ref{quad inv to ortho quad inv}
that one can calculate 
$B(\mu)^*\,a_{2r-1}$ and $B(\mu)^*\,d_T$
by determining the images of 
$w_{2r-1}= \pi^*(a_{2r-1})$ and of 
$\sum_{i\in T}w_{2i-1} v_{T-\{ i\}} 
= \pi^*(d_T)$ under $B(\mu)^*:H^*(BO(2n)) \to H^*(BO(2n))[w]$ 
and then expressing them in terms
of the generators 
(namely, $\pi^*(a_{2i-1})$, $\pi^*(d_S)$, $\pi^*(b_{4i})$, and 
$\theta(t)=w^2$) 
of the sub-algebra 
$\theta(H^*(BGO(2n))[t])\subset H^*(BO(2n))[w]$. 
This calculation is short in the case of the $a_{2r-1}$, 
and in the case of $d_T$ it can be done by a standard algorithm 
using Grobner bases. As an example, we have given 
the answer when $T$ has cardinality $2$.
\hfill$\square$ 

\medskip

\centerline{\bf The invariants $\lambda$ and $a_1$}

The elements $\lambda$ and $a_1$ of
$H^*(BGO(2n))$ are in $PH^*(BGO(2n))$ for all even ranks $2n$. 
The element $a_1$ is commonly known as the `discriminant', 
and is associated to the character
$\psi : GO(2n) \to \{ \pm 1\}$ defined by  
$g\mapsto \sigma(g)^n/\det(g)$
where recall that $\sigma : GO(2n) \to \C\,^*$ was defined by
the equality ${^t}g g = \sigma(g)I$.
It also has a well-known Gauss-Manin description in terms
of the fibration $V\to X$ where 
$V$ is the subvariety of $P(E)$ defined by the vanishing of $b$.
(So, $V$ is a quadric bundle in the original sense).

\centerline{\bf The ring of quadric invariants in rank $2$}

In rank $2$, we have $B(\mu)^*\lambda =  \lambda $, 
$B(\mu)^*a_1 = a_1$, and  
$B(\mu)^*\,b_4 = b_4 + (a_1^2 +\lambda) t + t^2$.
Hence the ring of quadric invariants is
$$PH^*(BGO(2)) = {\FF_2 [\lambda, a_1]\over (\lambda a_1)} \subset 
{\FF_2 [\lambda, a_1, b_4]\over (\lambda a_1)} = H^*(BGO(2))$$

\medskip

\centerline{\bf The ring of quadric invariants in rank $4$}

Generators for the subring 
$PH^*(BO(4)) \subset H^*(BO(4))$ of orthogonal quadric
invariants have been given by Toda in Proposition 3.12 of [T], which we
recall (correcting a minor misprint).  

\begin{proposition}\label{Toda rank 4} 
{\rm (Toda [T] Proposition 3.12)}
The elements $w_1$, $w_2^2 + w_1w_3$, and $w_1w_2w_3 + w_3^2 + w_1^2w_4$ 
generate the subring $PH^*(BO(4)) \subset H^*(BO(4))$.
\end{proposition}

To give a set of generators for the ring of quadric invariants 
$PH^*(BGO(4)) \subset H^*(BGO(4))$, 
we combine the above and Remark \ref{quad inv to ortho quad inv}, 
with the following table of the action of $B(\mu)^*$ in {\bf rank $4$}
which follows from Proposition \ref{action}.
\begin{eqnarray*} 
B(\mu)^*\lambda & = & \lambda \\ 
B(\mu)^*a_1     & = & a_1 \\  
B(\mu)^*a_3     & = & a_3 + a_1 t \\  
B(\mu)^*d_{1,2} & = & d_{1,2} + a_1 t^2  \\ 
B(\mu)^*\,b_4   & = & b_4 + a_1^2 t \\ 
B(\mu)^*\,b_8   & = & b_8 + (a_3^2 + \lambda^3 +  \lambda b_4 )t 
   + b_4 t^2 + a_1^2 t^3 + t^4  
\end{eqnarray*} 

\begin{proposition}  
The ring $PH^*(BGO(4))\subset H^*(BGO(4))$ of 
quadric invariants in rank $4$ 
is generated by the elements $\lambda$, $a_1$,  
$a_1a_3 + b_4$, and $a_3^2 + a_1d_{\{1, 2 \} }$. 
\end{proposition}  

\proof Let $P'$ be the subring of $H^*(BGO(4))$
generated by the elements $\lambda$, $a_1$, $a_1a_3 + b_4$, and 
$a_3^2 + a_1d_{\{1, 2 \} }$. We want to show that $P'=PH^*(BGO(4))$. 
From the above table of the effect of $B(\mu)^*$, it follows that
these elements are indeed quadric invariants, so 
$P'\subset PH^*(BGO(4))$. Further note that 
$\pi^*(\lambda) = 0 \in H^*(BO(4))$, while the images under $\pi^*$ of
$x=a_1$, $y=a_1a_3 + b_4$, and $z=a_3^2 + a_1d_{\{1, 2 \} }$ are 
respectively the generators $x'=w_1$, $y'=w_2^2 + w_1w_3$, and 
$z'=w_1w_2w_3 + w_3^2 + w_1^2w_4$ of the subring 
$PH^*(BO(4)) \subset H^*(BO(4))$.
Let $g \in H^*(BGO(4))$ be any quadric invariant. By the Remark
\ref{quad inv to ortho quad inv}, its image in 
$H^*(BO(4))$ is an orthogonal quadric invariant, so is expressible as a
polynomial $p(x',y',z')$ by Proposition \ref{Toda rank 4}. 
Hence $h = g - p(x,y,z)$ lies in the kernel of
$\pi^* : H^*(BGO(4)) \to H^*(BO(4))$. 

We claim that 
$$\ker(PH^*(BGO(4))\to PH^*(BO(4))) = \lambda k[\lambda, b_4]$$
If the claim is true then to prove the proposition we have to show that  
for any polynomial $p(\lambda,b_4)$, the element
$\lambda p(\lambda,b_4)$ lies in $P'$. This holds for 
the monomials of the form $\lambda ^i b_4^j$ 
for any $i>0$ and $j\ge 0$, in view of the equality
$\lambda ^i b_4^j=\lambda ^i y^j$, which holds as $\lambda a_1 =0$.
Hence the claim implies the proposition.

Now to prove the claim we consider any polynomial 
$g=g(\lambda,b_4,b_8)\in k[\lambda,b_4,b_8]$. We have to show that if it
lies in $\mbox{ker}(PH^*(BGO(4))\lra PH^*(BO(4)))$ then it lies in 
$\lambda k[\lambda, b_4]$.
We write the polynomial as $g=\sum_{i=0}^r p_i(\lambda ,b_4)b_8^i$, where 
$p_r(\lambda,b_4)\ne 0$. Consider the equalities
\begin{eqnarray*}
B(\mu)^*(\lambda g) &=&
B(\mu)^* \left(\sum_{i=0}^r (\lambda p_i(\lambda ,b_4))b_8^i \right) \\
&=& \sum_{i=0}^rB(\mu)^* (\lambda p_i(\lambda ,b_4))B(\mu)^*(b_8^i)  \\
&=& \sum_{i=0}^r (\lambda p_i(\lambda ,b_4))(B(\mu)^*(b_8))^i.
\end{eqnarray*}
Using the fact that  
$B(\mu)^*(b_8) = 
b_8+a_1^2t^3+b_4t^2+(a_3^2+\lambda ^3+\lambda b_4)t^3+t^4$,   
we find that the coefficient of $t^{4r}$ in the expression for
$B(\mu)^*(\lambda g)$ is exactly equal to $\lambda p_r(\lambda,b_4)$,
which is a contradiction unless $r=0$.
Hence $g=p_0(\lambda,b_4)$. 
This completes the proof of the claim, hence the proposition.
\hfill$\square$ 

\medskip

\bigskip

\centerline{\bf The ring $PH^*(BGO(4m+2))$ 
of quadric invariants in rank $4m+2$}

We first recall the following short exact sequence which is used 
in [H-N] to compute the cohomology ring of $BGO(2n)$ 
(in the present case, $BGO(4m+2)$)
$$ 0 \to \lambda \, \FF_2[b_4,\ldots,b_{8m+4}] 
\to H^*(BGO(4m+2)) \to C \to 0 $$ 
where $C$ is the subring of $H^*(BO(4m+2))$ 
which is the image of $\pi^* : H^*(BGO(4m+2)) \to H^*(BO(4m+2))$.
We recall from [H-N] that $C$ equals the kernel of the derivation 
$d_1 = s = \sum w_{2i-1}{\pa \over \pa w_{2i}}$
on $H^*(BO(4m+2)) = \FF_2[w_1,\ldots,w_{4m+2}]$.

Recall that by Remark \ref{quad inv to ortho quad inv}, 
for $A = H^*(BO(4m+2))$ and $A' = H^*(BGO(4m+2))$,
the map $\pi^* : A'\to A$ maps $PA' \to PA$.
Also recall from [H-N] that the map $\pi^* : A'\to A$ is 
injective in odd ranks.

\medskip

\centerline{\bf Generators ${\alpha'}_{2i-1}$, ${\delta'}_T$ 
and ${\beta'}_{4i}$ for the primitive ring $PH^*(BGO(4m+2))$}

Consider the maps $d_j: A\to A$ for $j \ge 1$,
defined by $d_j(w_r) = {4m+2 \choose j}w_{r-j}$.  
By definition of the $d_j$, we have $PA \subset \ker(d_j)$ for
$j \ge 1$, in particular, 
the elements ${\alpha}_{2i-1}$ and
${\delta }_T \in PH^*(BO(4m+2))$ are in the kernel 
of the map $d_j$ for $j \ge 1$. 
Hence these elements lie in $C = \ker(d_1) = \im(\pi^*)$.
As the map $\pi^*$ is injective in odd ranks, 
there exist unique elements
$${\alpha'}_{2i-1} \in  H^{2i-1}(BGO(4m+2)) \mbox{ for } 1\le i \le 2m+1$$ 
and 
$${\delta'}_T \in  H^{2\sum p_i -1}(BGO(4m+2)) 
\mbox{ for } T=\{p_1,\ldots,p_r\} \subset \{2,\ldots,2m+1 \}, ~ r\ge 2$$ 
such that they map to ${\alpha}_{2i-1}$ 
and ${\delta }_T$ respectively.  
By the injectivity of $\pi^*$ in odd dimensions, 
these elements actually lie in $PH^*(BGO(4m+2))$.

Recall the definition of the elements $\beta _{4i} \in PA$ by 
$$\beta _{4i}
= {\wh w}_{2i}^2+{\wh w}_{2i-1}({\wh w}_{2i}w_1+{\wh w}_{2i-1}w_2)
~~\mbox{for } 2 \le i \le 2m+1$$ 
As ${\wh w}_{2i}w_1+{\wh w}_{2i-1}w_2 = d_1(\wh{w}_{2i}w_2)$
and $d_1\circ d_1 =0$, 
${\wh w}_{2i}w_1+{\wh w}_{2i-1}w_2 \in \ker(d_1) = \im(\pi^*)$.  
As ${\wh w}_{2i}w_1+{\wh w}_{2i-1}w_2 \in H^{2i+1}(BO(4m+2))$
is of odd rank $2i+1$, there exist a unique element 
$$g_{2i+1} \in H^*(BGO(4m+2)) \mbox{ with }
\pi^*(g_{2i+1}) = {\wh w}_{2i}w_1+{\wh w}_{2i-1}w_2$$
Recall from the Appendix the definition of elements
${\wh c}_i \in H^*(BGL(4m+2))$ for $i>2$,
following Toda [T]. 
We now define elements ${\wh b}_{4i}$ by 
$${\wh b}_{4i} = p^*({\wh c}_{2i}) \in H^{4i}(BGO(4m+2))
~~\mbox{for } 2 \le i \le 2m+1$$ 
where 
$p^*:H^*(BGL(4m+2)\to H^*(BGO(4m+2))$ is the map induced by the 
inclusion $GO(4m+2)\hra GL(4m+2)$.
Then we observe that $\pi^*({\wh b}_{4i})={\wh w}_{2i}^2$,
as $\wh{c}_i \mapsto \wh{w}_i^2$ as already seen.
We define the element ${\beta'}_{4i} \in H^{4i}(BGO(4m+2))$ 
by the equality
$${\beta'}_{4i}={\wh b}_{4i}+{\alpha'}_{2i-1}g_{2i+1}
~~\mbox{for } 2 \le i \le 2m+1$$
so that we have
$$\pi^*({\beta'}_{4i})=\beta _{4i}$$

\begin{lemma}\label{surjectivity} 
The elements ${\beta'}_{4i}$ lie in 
$PH^*(BGO(4m+2)) \subset H^*(BGO(4m+2))$.
\end{lemma}

\proof By definition of 
$g_{2i+1}$, we have $\pi^*(g_{2i+1}) = {\wh w}_{2i}w_1+{\wh w}_{2i-1}w_2$.
Hence,
\begin{eqnarray*}
B(\mu)^*(\pi^*(g_{2i+1}))
& = & B(\mu)^*({\wh w}_{2i}w_1+{\wh w}_{2i-1}w_2)\\
& = &({\wh w}_{2i}w_1+{\wh w}_{2i-1}w_2)+w^2{\wh w}_{2i-1}
\end{eqnarray*}

Therefore by the commutativity of the square
$$\begin{array}{ccc}
H^*(BGO(2n)) & \stackrel{B(\mu)^*}{\lra} & H^*(BGO(2n))[t]\\
\dna & & \dna \\
H^*(BO(2n)) & \stackrel{B(\mu)^*}{\lra} & H^*(BO(2n))[w]
\end{array}$$
we see that
$B(\mu)^*(g_{2i+1})=g_{2i+1}+t{\alpha'}_{2i-1}+t \lambda f$
for some $f \in H^{2i - 3}(BGO(2n))[t]$.
As $\lambda {\alpha'}_{2i-1} = 0$, 
$B(\mu)^*({\alpha'}_{2i-1}g_{2i+1}) =
{\alpha'}_{2i-1}g_{2i+1}+t{\alpha'}_{2i-1}^2$.
This equality is used in the third line of the following chain of 
equalities:

\begin{eqnarray*}
B(\mu)^*({\beta'}_{4i})
& = & B(\mu)^*({\wh b}_{4i} + {\alpha'}_{2i-1}g_{2i+1}) 
~\mbox{ by definition of } {\beta'}_{4i},\\
& = & B(\mu)^*( p^*({\wh c}_{2i}))+B(\mu)^*({\alpha'}_{2i-1}g_{2i+1}) 
~\mbox{ by definition of } {\wh b}_{4i},\\
& = & p^*(B(\mu)^*({\wh c}_{2i}))+{\alpha'}_{2i-1}g_{2i+1}
      + t{\alpha'}_{2i-1}^2 ~\mbox{ as explained above}, \\
& = & {\wh b}_{4i}+{\alpha'}_{2i-1}g_{2i+1}+t(p^*({\wh c}_{2i-1})
      + {\alpha'}_{2i-1}^2)
\end{eqnarray*}
The last line above follows from the fact that 
$B(\mu)^*({\wh c}_{2i})={\wh c}_{2i}+t{\wh c}_{2i-1}$ and that
$p^* \circ B(\mu)^*=B(\mu)^* \circ p^*$.

Hence to complete the proof of the Lemma \ref{surjectivity},
it is sufficient to prove the following. 

\begin{lemma}\label{cbar=b2} 
With the above notations, we have 
$p^*({\wh c}_{2i-1})={\alpha'}_{2i-1}^2$ for $2 \le i \le 2m+1$.
\end{lemma}

{\bf Proof} of Lemma \ref{cbar=b2}: We first observe that 
$B(\mu)^*(p^*({\wh c}_{2i-1}))
= B(\mu)^*({\alpha'}_{2i-1}^2)={\wh w}_{2i-1}^2$.
Hence $p^*({\wh c}_{2i-1})={\alpha'}_{2i-1}^2+\lambda h$,
where $h = h(\lambda,b_4,\ldots,b_{8m+4})$.
We just have to show that $h = 0$.

For that we consider the group homomorphisms
$\tau : O(2) \times O(2m+1) \to O(4m +2)$,
$\tau' : GO(2) \times O(2m+1) \to GO(4m +2)$ and 
$\tau'' : GL(2) \times GL(2m+1) \to GL(4m +2)$, 
induced by the tensor product 
$\C^2\otimes\C^{2m+1}\to \C^{4m+2}$. 
These fit in the 
following commutative diagram of group homomorphisms, 
in which the vertical maps are natural
inclusions.

$$\begin{array}{ccc}
O(2) \times O(2m+1)  &   \stackrel{\tau}{\to} & O(4m +2) \\
\dna                 &                        & \dna      \\
GO(2) \times O(2m+1) &  \stackrel{\tau'}{\to} & GO(4m +2) \\
\dna                 &                        & \dna      \\
GL(2) \times GL(2m+1)& \stackrel{\tau''}{\to} & GL(4m +2) 
\end{array}$$

Hence we get the following commutative diagram of ring homomorphisms.

$$\begin{array}{ccc}
H^*(BGL(4m+2))&\stackrel{B({\tau''})^*}{\to} & 
H^*(BGL(2))\otimes H^*(BGL(2m+1))\\
{\scriptstyle p^*} \downarrow & & \downarrow {\scriptstyle 
p^* \otimes (p \pi)^*} \\
H^*(BGO(4m+2))&\stackrel{B({\tau'})^*}{\to} & 
H^*(BGO(2))\otimes H^*(BO(2m+1))\\
{\scriptstyle \pi^*} \downarrow & & \downarrow 
{\scriptstyle \pi^* \otimes \id }~\\
H^*(BO(4m+2))&\stackrel{B(\tau)^*}{\to}& H^*(BO(2))\otimes H^*(BO(2m+1))
\end{array}$$

By Propositions \ref{tau=0} and \ref{tau''=0} we have
$B({\tau''})^*({\wh c_{2i-1}})=0$ and 
$B({\tau'})^*(\alpha'_{2i-1})=0$.
These facts imply that $B({\tau''})^*(\lambda h)=0$.
From this, we can deduce that $\lambda h=0$ by means of the following lemma
(statement \ref{inj}.(3) below),
completing the proof of the Lemma \ref{cbar=b2}.

\begin{lemma}\label{inj} 
With the above notations we have

(i) $B({\tau'})^*(\lambda)= \lambda \otimes 1$,
where we also denote by $\lambda$ the class $\ov{c}_1(L)\in H^*(BGO(2))$, 

(ii) The map $B({\tau'})^*$ is injective on 
$\FF_2[\lambda, b_4,\ldots,b_{8m+4}]$.
\end{lemma}

{\bf Proof} of Lemma \ref{inj}:  
To prove the first part, we use the commutativity of the square  
$$\begin{array}{clc}
H^*(BGO(4m+2))& \stackrel{B({\tau'})^*}{\to} & 
                                        H^*(BGO(2))\otimes H^*(BO(2m+1))\\
{\scriptstyle \pi^*}\downarrow&&\downarrow{\scriptstyle \pi^*\otimes\id}\\
H^*(BO(4m+2))&\stackrel{B(\tau)^*}{\to}& H^*(BO(2))\otimes H^*(BO(2m+1))
\end{array}
$$
and the fact that the map 
{\small
$$
\oplus_{i=0}^1(H^i(BGO(2))\otimes H^{2-i}(BO(2m+1)))
\to \oplus_{i=0}^1(H^i(BO(2))\otimes H^{2-i}(BO(2m+1)))
$$
}
is injective. This implies that the image of $\lambda$ under the map
$B({\tau'})^*$ is contained in 
$H^2(BGO(2))\otimes H^0(BO(2m+1))$. 
But the K\"unneth projection 
$$H^*(BGO(2))\otimes H^*(BO(2m+1))\to H^*(BGO(2))$$
composed with $B({\tau'})^*$ is a map 
$H^*(BGO(4m+2))\to H^*(BGO(2))$ which is induced by the inclusion
$GO(2) \subset GO(4m+2)$ defined by 
$$A \mapsto 
\left( \begin{array}{lcl}
A &       &  \\
  &\ddots &  \\  
  &       & A
\end{array} \right)$$
Hence the map $H^*(BGO(4m+2))\to H^*(BGO(2))$ takes 
$\lambda \mapsto \lambda$. This implies that 
$B({\tau'})^*(\lambda )= \lambda \otimes 1$.
This completes the proof of the first assertion.

For proving statement (ii), first
consider the commutative square 
$$\begin{array}{clc}
H^*(BGO(4m+2))&\stackrel{B({\tau'})^*}{\to} & 
H^*(BGO(2))\otimes H^*(BO(2m+1))\\
{\scriptstyle \pi^*} \downarrow & & \downarrow {\scriptstyle B(u)^*} \\
H^*(BO(4m+2))&\stackrel{B(i)^*}{\to}& H^*(BO(2m+1))
\end{array}$$
where $B(i)^*$ is induced by the the inclusion
$$i : O(2m+1)\hra O(4m+2) : A \mapsto \left( 
\begin{array}{ll} 
  A &   \\  
    & A 
\end{array}
\right)$$ 
and $u : O(2m+1) \to GO(2) \times O(2m+1)$ is the homomorphism
defined by $g\mapsto (1,g)$. 
It follows from the definition of $B(i)^*$ that
$B(i)^*(w_{2j})=w_j^2$. 
Hence for any polynomial $p(w_2,\ldots,w_{4m+2})$  we have 
$B(i)^*(p)=p(w_1^2,\ldots,w_{2m+1}^2)$.
By the commutativity of the above square it follows that for 
any polynomial $p=p(b_4,\ldots,b_{8m+4})\in H^r(BGO(4m+2))$
we have $B({\tau'})(p)=1\otimes p(w_1^4,\ldots,w_{2m+1}^4)+q$
where $q \in \oplus _{j=1}^rH^j(BGO(2))\otimes H^{r-j}(BO(2m+1))$.

Now consider a polynomial 
$p = p(\lambda,b_4,\ldots,b_{8m+4}) \in H^{2r}(BGO(4m+2))$. 
Writing it as $p = \sum_{i=0}^r \lambda ^i p_i(b_4,\ldots,b_{8m+4})$,
let $\ell$ be the smallest integer such that \\ 
$p_{\ell}(b_4,\ldots,b_{8m+4}) \ne 0$.
Then we see that 
$B({\tau'})(p)=\lambda ^{\ell} \otimes p_{\ell}(w_1^4,\ldots,w_{2m+1}^4)+q$,
where $q \in \oplus _{j=\ell+1}^{2r}H^j(BGO(2))\otimes H^{2r-j}(BO(2m+1))$.
Hence $B({\tau'})(p) \ne 0$.
This proves Lemma \ref{inj}.(ii), and hence completes the
proofs of the Lemmas \ref{cbar=b2} and \ref{surjectivity}. 
\hfill $\square$

\begin{proposition} 
The map $PH^*(BGO(4m+2))\to PH^*(BO(4m+2))$ is surjective.
\end{proposition}

\proof We know that $PH^*(BO(4m+2))$ is generated by elements of the form
$\alpha_{2i-1}$, $\beta_{4i}$ and $\delta_T$. These are the images of
$\alpha'_{2i-1}$, $\beta'_{4i}$ and $\delta'_T$ in $PH^*(BGO(4m+2))$,
so the proposition follows.
\hfill $\square$

\begin{theorem}\label{maingo} 
The primitive subring $PH^*(BGO(4m+2)) \subset H^*(BGO(4m+2))$ 
is generated by the finite set consisting of the following elements. 

(i) The element $\lambda$,

(ii) the element $\alpha'_1 = a_1$, 

(iii) the elements 
${\alpha'} _{2i-1}$ for $2 \le i \le 2m+1$,

(iv) the elements ${\beta'} _{4j}$ for $2 \le j \le 2m+1$,
and

(v) the elements ${\delta'}_T$, for
$T \subset \{2,\ldots,2m+1 \}$ with cardinality $|T| \ge 2$.
\end{theorem}

\proof Let $P'$ be the subring of $PH^*(BGO(4m+2))$ generated by the 
elements listed in the statement of the theorem. 
We want to show that $P' = PH^*(BGO(4m+2))$.
By Proposition \ref{surjectivity} the map 
$PB(\mu)^* : PH^*(BGO(4m+2)) \to PH^*(BO(4m+2))$
is surjective. As its kernel is 
$\lambda \FF_2[\lambda, b_4, \ldots, b_{8m+4}]\cap PH^*(BGO(4m+2))$,
it is enough to show that we have the following inclusion of sets: 
$$\lambda \FF_2[\lambda, b_4, \ldots, b_{8m+4}]\cap PH^*(BGO(4m+2)) 
\subset P'$$
Consider any polynomial $p=p(\lambda, b_4,\ldots,b_{8m+4})$.
We can write $b_{4i}$ as ${\wh b}_{4i}+z$ where $z$ is a polynomial in 
elements of $H^j(BGO(4m+2))$ for $j<4i$ (this holds because
the elements $b_{4i}$ and ${\wh b}_{4i}$ come from the elements
$\ov{c}_{2i}$ and $\wh{c}_{2i}$ of $H^{4i}(BGL(4m+2))$,
and the corresponding statement holds in $H^*(BGL(4m+2))$).

As $\lambda$ annihilates 
odd degree elements of $H^*(BGO(2n))$, we can assume that
$z$ is a polynomial in the elements $\lambda, b_4,\ldots,b_{4i-4}$.
We substitute this in $p$ to conclude, by iteration, 
that the polynomial $\lambda p$
is an element of 
$\lambda \FF_2[\lambda, b_4, {\wh b}_8,\ldots,{\wh b}_{8m+4}]$.
It follows from the definition of ${\beta'}_{4i}$ that
$\lambda {\beta'}_{4i}=\lambda {\wh b}_{4i}$. 
This implies that for any polynomial 
$p=p(\lambda, b_4, \ldots, b_{8m+4})$, we have
$\lambda p \in 
\lambda \FF_2[\lambda, b_4, \beta'_8,\ldots,\beta'_{8m+4}]$.
Now let 
$p = \sum_{i=0}^r b_4^i p_i(\beta'_8,\ldots,\beta'_{8m+4}))
\in H^*(BGO(4m+2))$
be such that 
$\lambda p \in 
\lambda \FF_2[\lambda, b_4, \ldots, b_{8m+4}]\cap PH^*(BGO(4m+2))$.
Suppose $r \ge 1$.
As $B(\mu)^*(p_i)=p_i$ and $B(\mu)^*(b_4)=b_4+a_1^2t+t^2$,
the coefficient of $t^{2r}$ in the expansion of 
$B(\mu)^*(\sum_{i=0}^r (\lambda p_i) b_{4}^i)$
is equal to $\lambda p_r(\beta'_8,\ldots,\beta'_{8m+4})$. 
Hence if $\lambda p$ is in $PH^*(BGO(4m+2))$, we must have
$\lambda p_r(\beta'_8,\ldots,\beta'_{8m+4})=0$. 
Continuing this way, we get 
$\lambda p = \lambda \cdot p_0(\beta'_8,\ldots,\beta'_{8m+4})$.
This lies in $P'$, which shows the desired inclusion and 
completes the proof of the theorem.  \hfill $\square$

\rem 
In fact the method of the proof gives us the 
short exact sequence 
$$ 0\to \FF_2[\lambda,{\beta'}_8,\ldots,{\beta'}_{8m+4}]
\stackrel{\lambda}{\to}
PH^*(BGO(4m+2))\to PH^*(BO(4m+2))\to 0 $$
which can be used to write down the relations in the ring
$PH^*(BGO(4m+2))$. Similarly, we can write the relations between our 
generators of $PH^*(BO(4))$. 
We do not include the description of relations
in this paper, as it does not contribute here to our main theorems 
on degenerating quadric bundles.

%

\section{The Gysin boundary map}

\begin{lemma}\label{cupgysin} 
Let $Y$ be a  Hausdorff topological space, and 
let $\pi : N\to Y$ be a real vector bundle bundle on $Y$,
of real rank $r$. Let $N_o \subset N$ denote the complement of 
the zero section
$i : Y\hra N$ of $N\to Y$.  
With the above notation, the following diagram is commutative
$$\begin{array}{ccccc} 
H^p(N) & \times & H^q(N_o) & \stackrel{\cup}{\lra} & H^{p+q}(N_o) \\
{\scriptstyle i^*}\dna &&{\scriptstyle\pa}\dna&&{\scriptstyle\pa}\dna\\
H^p(Y) & \times & H^{q-r+1}(Y) & \stackrel{\cup}{\lra} & H^{p+q-r+1}(Y) 
\end{array}$$
In other words, given $a \in H^p(N)$ and $b\in H^q(N_o)$, we have
$\pa(a \cup b) = i^*(a) \cup \pa(b) \in H^{p+q-r+1}(Y)$.
\end{lemma}

\proof Let $\t\in H^r(N,N_o)$ be the Thom class, 
so that for each $p$ we have the Thom isomorphism $T$, which is the 
composite 
$$T: H^p(Y)\stackrel{\pi^*}{\lra}  H^p(N)
\stackrel{-\cup \t }{\lra} H^{p+r}(N,N_o)$$ 
where both $\pi^*$ and $-\cup\tau$ are isomorphisms.
The Gysin boundary map $\pa: H^p(N_o) \to  H^{p-r+1}(Y)$
is the composite
$$\pa: H^p(N_o) \stackrel{\delta}{\to}
H^{p+1}(N,N_o)\stackrel{T^{-1}}{\to}
H^{p-r+1}(Y)$$
of the connecting homomorphism $\delta$ for the pair $(N,N_o)$ 
with the inverse of
the Thom isomorphism $T$.
If $\alpha\in H^p(N)$ and $\beta\in H^{q-r+1}(N)$,
then by associativity of cup product,
$(\alpha\cup\beta)\cup\t = \alpha\cup(\beta\cup\t)$, in other
words, the following diagram is commutative.
$$\begin{array}{ccccc} 
H^p(N)&\times&H^{q-r+1}(N)&\stackrel{\cup}{\lra}&H^{p+q-r+1}(N)\\ 
{\scriptstyle\id}\dna &&{\scriptstyle\cup\t}\dna&&{\scriptstyle\cup\t}\dna\\
H^p(N) & \times & H^{q+1}(N,N_o) & \stackrel{\cup}{\lra} & 
H^{p+q+1}(N,N_o)
\end{array}$$
By composing with the 
isomorphism $\pi^* :H^*(Y)\to H^*(N)$, this gives 
the following commutative diagram, where $T$ is the 
Thom isomorphism.
$$\begin{array}{ccccc} 
H^p(Y)&\times&H^{q-r+1}(Y)&\stackrel{\cup}{\lra}&H^{p+q-r+1}(Y)\\ 
{\scriptstyle \pi^*}\dna &&{\scriptstyle T}\dna&&{\scriptstyle T}\dna\\
H^p(N) & \times & H^{q+1}(N,N_o) & \stackrel{\cup}{\lra} & 
H^{p+q+1}(N,N_o)
\end{array}$$
Hence as $i^*$ is the inverse of $\pi^*$,
the following diagram commutes.
$$\begin{array}{ccccc}
H^p(N) & \times & H^{q+1}(N,N_o) & \stackrel{\cup}{\lra} & 
H^{p+q+1}(N,N_o)\\
{\scriptstyle i^*}\dna &&{\scriptstyle T^{-1}}\dna&
                                      &{\scriptstyle T^{-1}}\dna\\
H^p(Y)&\times&H^{q-r+1}(Y)&\stackrel{\cup}{\lra}&H^{p+q-r+1}(Y)
\end{array}$$
For any integer $p$, let 
$\delta :H^p(N_o) \to H^{p+1}(N,N_o)$ denote the connecting homomorphism
for the singular cohomology of the pair $(N,N_o)$.
For any integers $p$ and $q$, consider the cup product maps
$H^p(N) \otimes H^q(N_o) \to H^{p+q}(N_o)$
and $H^p(N)\otimes H^{q+1}(N,N_o) \to H^{p+q+1}(N,N_o)$.
It is known that these fit in a commutative diagram
$$\begin{array}{ccccc} 
H^p(N) & \times & H^q(N_o) & \stackrel{\cup}{\lra} & H^{p+q}(N_o) \\
{\scriptstyle\id}\dna &&{\scriptstyle\delta}\dna&&{\scriptstyle\delta}\dna\\
H^p(N) & \times & H^{q+1}(N,N_o) & \stackrel{\cup}{\lra} & 
H^{p+q+1}(N,N_o) 
\end{array}$$
Now the lemma follows by juxtaposing the above two commutative diagrams.
\hfill $\square$

\rem 
The above lemma is for coefficients $\FF_2$. A version
with arbitrary coefficient ring $R$ is possible, where we have to assume that
$N$ is orientable with respect to coefficients $R$, 
and the conclusion is that the the diagram in statement \ref{cupgysin} is 
graded-commutative, that is, 
$\pa(a \cup b) = (-1)^{\deg(a)} i^*(a) \cup \pa(b)$.

\rem\label{remcupGysin} 
We will apply the above lemma in the following situation. 
Let $X$ be a topological space and $i:Y\hra X$ a closed subspace which is
a topological divisor in $X$, that is, 
there exists a rank $2$ real vector bundle $\pi: N\to Y$ over $Y$,
together with a continuous map $\varphi: N\to X$
which is a homeomorphism of $N$ with an open neighbourhood
$U$ of $Y$ in $X$, such that $\varphi$ restricted to the 
zero section $Y\subset N$ is $\id_Y$.
Composing the restriction $ H^p(X-Y) \to H^p(U-Y)$ with the map  
$H^p(U-Y)\stackrel{\varphi^*}{\to} H^p(N_o)$ and 
the Gysin boundary map
$\pa : H^p(N_o) \to H^{p-1}(Y)$, we get a map
$$\pa : H^p(X-Y) \to H^{p-1}(Y)$$
which is by definition the Gysin boundary map for the pair $(X,Y)$.
It can be shown that the map $\pa$ is independent of the
choice of $\pi: N\to Y$ and $\varphi$. 
The Lemma \ref{cupgysin} implies that 
given $a \in H^p(X)$ and $b\in H^q(X-Y)$, we have
$$\pa(a \cup b) = i^*(a) \cup \pa(b) \in H^{p+q-r+1}(Y)$$

\rem\label{gysin.summary} 
We now summarize the properties of the Gysin boundary map that we will use.
Let $L\to Y$ be a complex line bundle.
Let $L_o\subset L$ be the complement of the zero section, and
let $\pi : L_o \to Y$ be the projection.
Let $\pa : H^*(L_o)\to H^{*-1}(Y)$ be the Gysin boundary,
and $\pi^* : H^*(Y) \to H^*(L_o)$ 
the pullback under $\pi : L_o \to Y$.
Let $s = \pi^*\circ \pa: H^*(L_o) \to H^{*-1}(L_o)$. 
Let $\mu : \C\,^*\times L_o \to L_o$ be the scalar multiplication, and let
$p : \C\,^*\times L_o \to L_o$ be the projection.
Let $\eta\in H^1(\C\,^*)$ be the generator.
Let $\sq^i$ denote the $i$ th Steenrod operation.
Then for any $x,\,x'\in H^*(L_o)$ 
and $y\in H^*(Y)$, we have the following 
basic equalities. 
\begin{eqnarray}
\pa(x\cup \pi^*y) & = & \pa(x)\cup y\\
(\mu^* -p^*)(x) & = & \eta\otimes s(x)\\
\pa\sq^ix & = &  \sq^i \pa x \\ 
\pa(x^2) & = & 0 \\
s(x\cup x') & = & s(x)\cup x'+ x\cup s(x')
\end{eqnarray}
The property (1) is given by Lemma \ref{cupgysin}.
The properties (2)-(5) are proved in Section 2 of [H-N].
Note that property (5) can be expressed by saying that
the map $s$ is a derivation on the cohomology ring
$H^*(L_o)$. The above equalities are used only in the explicit
computations of the Gysin boundary images of the quadric invariants
(see Section 8).

%
%
%

\section{Topological vanishing multiplicity}

\centerline{\bf Definition of vanishing multiplicities}

Let $N$ be a $2$-dimensional real vector space, together with an orientation,
and let $L$ be a $1$-dimensional complex vector space.
Let $f : N-\{ 0\} \to L -\{ 0\}$ be a map. 
We denote by $\nu_0(\sigma) \in \Z$ the winding number of this map. 
Its sign depends on the orientation on $N$,
but its absolute value $|\nu_0|(f)\in \Z^{\ge 0}$ and its parity
$\ovnu_Y(f) \in \{ 0,1\} \subset \Z^{\ge 0}$
are independent of the orientation chosen.
Here the parity $\ovnu_Y(f)$ is defined to be $0$
or $1$ depending on whether $|\nu_0|(f)$ is even or odd.

Next, we globalize the above. Let $Y$ be a Hausdorff space,
of the homotopy type of a CW complex,  
and let $\pi :N\to Y$ be a real rank $2$ vector bundle on $Y$. 
Let $L$ be a complex
line bundle on the total space of $N$. Let $N-Y$ be the complement of the 
zero section $Y\subset N$, and let $\sigma \in \Gamma(N-Y,L)$ be a section, 
which is nowhere vanishing. Choose an isomorphism 
$\theta : L \to \pi^*(L_Y)$ of complex line bundles, 
that is identity over $Y$ (which exists by Remark \ref{(2)}). 
For any $y\in Y$,
$\sigma$ gives a map $f_y : N_y - \{ 0\} \to L_y - \{ 0\}$. 
The numbers $|\nu_0|(f_y)\in \Z^{\ge 0}$ and $\ovnu_0(f_y) \in 
\{ 0,\,1 \}\subset \Z$ are well-defined, independent of the choice of 
the isomorphism $\theta : L \to \pi^*(L_Y)$ as follows from
Remark \ref{(2)}. Hence we regard them as functions
on $Y$. As such, these are constant on
path components of $Y$, so give functions
$$|\nu_Y|(\sigma) : \pi_0(Y) \to \Z^{\ge 0} 
\mbox{ and } \ovnu_Y(\sigma) : \pi_0(Y) \to \{0,\,1\} \subset \Z$$
If $N$ is orientable, then a choice of orientation gives rise to a function
$$\nu_Y(\sigma) : \pi_0(Y) \to \Z$$
similarly defined. 

\defin\label{(3)} 
By definition of $0$ th singular cohomology, we thus have an
element 
$|\nu_Y|(\sigma) \in H^0(Y;\Z)$ which we call the 
{\bf absolute topological vanishing multiplicity}, an element
$\ovnu_Y(\sigma) \in H^0(Y;\Z)$ which we call the 
{\bf topological vanishing parity}, and when $N$ is
given an orientation, an element
$\nu_Y(\sigma) \in H^0(Y;\Z)$ which we call the 
{\bf oriented topological vanishing multiplicity} of the non-vanishing
section $\sigma \in \Gamma(N-Y, L)$. 

\rem\label{tubular neighbourhood} 
More generally, if $X$ is a 
paracompact, Hausdorff topological space, $Y\subset X$
is a topological divisor, $L$ a complex line bundle on $X$, and
$\sigma \in \Gamma(X-Y,L)$ is a nowhere vanishing section, 
we can define the 
elements $|\nu_Y|(\sigma) \in H^0(Y;\Z)$, 
$\ovnu_Y(\sigma) \in H^0(Y;\Z)$
by choosing a pair $(N,\varphi)$ where $N$ is a 
rank $2$ real vector bundle on 
$Y$ and $\varphi : N \to X$ is an open embedding which 
maps the zero section 
of $N$ identically to $Y$, and then applying the above definition to
the pullback under $\varphi$. The resulting
topological vanishing multiplicity $|\nu|$ and its parity $\ovnu$ 
are well-defined, independent of the choice of 
the tubular neighbourhood $(N,\varphi)$. Moreover, if the 
bundle $N$ is orientable, then the choice of an
orientation gives $\nu_Y(\sigma) \in H^0(Y;\Z)$.

For a definition of $|\nu_Y|(\sigma)$ in terms of 
local cohomology, see [N].

\rem\label{(2)} 
As the zero section $Y\subset N$ 
is a strong deformation retract of $N$, 
for any complex line bundle $L$ on the total space $N$
there exists an isomorphism of complex line bundles
$\theta : L \to \pi^*(L_Y)$ which is identity on the zero section
$Y\subset N$, where $L_Y= L|_Y$, and 
any two such isomorphisms 
$\theta_1,\theta_2 : L \to \pi^*(L_Y)$ 
are related by $\theta_1 =f\theta_2$ 
where $f : N \to \C\,^*$ is a function such that 
$f$ is homotopic to $1$ relative to $Y$
(in symbols, $f \sim 1$ rel $Y$).

\rem\label{(4)} 
Let $\pi : N \to Y$ will be a real vector bundle of rank $2$ on a space $Y$,
equipped with a fiber-wise Riemannian metric (that is, a continuous
positive definite symmetric real bilinear form on $N$).
Let $\sigma\in \Gamma(N-Y,L)$
be a nowhere vanishing section with $|\nu_Y|(\sigma) \ne 0$.
Then $N$ is orientable, and has a unique orientation such that 
$\nu_Y(\sigma)$ is positive. This orientation, together with the given
metric, defines a $\C$-linear structure on $N$. 
Using this, we can regard $N$ as a complex line bundle.

\begin{lemma}\label{BZ} 
Let $N\to Y$ be a rank $2$ real vector bundle on $Y$, let
$D\subset N$ be the disk bundle over $Y$ of radius $1$ with 
respect to a  metric on 
the vector bundle $N$, and let $\pa D$ denote the boundary of $D$.
Let $g : \pa D \to \C\,^*$ be a continuous function
such that the restriction of $g$ to any fiber of
$p: \pa D\to Y$ is homotopic to a constant function. 
Then $g : \pa D \to \C\,^*$ 
admits a prolongation $\ov{g} : D \to \C\,^*$.
\end{lemma}

\proof By working on one component of $Y$ at a time, we can assume that
$Y$ is connected, that is, $\pi_0(Y)=1$.
Consider the homotopy exact sequence 
$$\pi_1(S^1) \stackrel{i_*}{\to} \pi_1(\pa D) 
\stackrel{p_*}{\to} \pi_1(Y) \to 1$$
for the fibration $\pa D \to Y$, where $i: S^1\hra \pa D$ 
denotes one of the fibers. The hypothesis of the lemma shows that 
the homomorphism $g_* : \pi_1(\pa D) \to \pi_1(\C\,^*)$ 
is trivial on the image of
$\pi_1(S^1) \to \pi_1(\pa D)$, so factors through $\pi_1(\pa D) 
\stackrel{p_*}{\to} \pi_1(Y)$ giving $\theta : \pi_1(Y) \to \pi_1(\C\,^*)$.
By identifying $\pi_1(D)$ with $\pi_1(Y)$ via the projection $D\to Y$,
we have a group homomorphism $\theta : \pi_1(D) \to \pi_1(\C\,^*)$.
As $\C\,^*$ is a $K(\pi, 1)$-space, every 
group homomorphism $\pi_1(X)\to \pi_1(\C\,^*)$ for a CW complex $X$
is of the form $f_*$ for some continuous function $f :X\to \C\,^*$ which is 
unique up to homotopy. We can therefore choose a function
$f : D \to \C\,^*$ with $f_* = \theta : \pi_1(D) \to \pi_1(\C\,^*)$.
On $\pa D$, we have $g_* = (f|\pa D)_*$, so by the above there exists
a homotopy $H: \pa D \times I\to \C\,^*$ from $g$ to $f|\pa D$.
Now let $D_{1/2}\subset D$ denote the disk bundle of radius half,
and let $D^o_{1/2}$ denote its interior. Note that we have
a canonical homeomorphism $\psi : \pa D \times I\to D - D^o_{1/2}$,
defined by the formula $\psi(v,t) = (1 - t/2)v$. 
Now define $\ov{g} : D \to \C\,^*$ by putting
$$\ov{g}(v) = \left\{ 
\begin{array}{ll} 
H\circ \psi^{-1}(v) & ~~\mbox{for } |\!|v|\!| \ge 1/2, \\ 
f                   & ~~\mbox{for } |\!|v|\!| \le 1/2. 
\end{array}\right. 
$$
This proves the lemma. \hfill$\square$

\begin{lemma}\label{non-vanishing prolongation} 
{\bf (Non-Vanishing prolongation) } 
Let $N\to Y$ be a rank $2$ real vector bundle on $Y$
equipped with a metric, 
let $L$ be a complex line bundle on $N$, and let
$\sigma \in \Gamma(N,L)$ be a global section, which 
is non-vanishing outside $Y$, with $|\nu_Y|(\sigma)=0$.  
Then the line bundle $L$ is trivial on $N$, and moreover
there exists a non-vanishing global section 
$\sigma'\in \Gamma(N,L)$ such that $\sigma'$ and $\sigma$ coincide
on the open set $N-D$, where 
$D\subset N$ be the disk bundle over $Y$ of radius $1$.
\end{lemma} 

\proof We first prove that $L$ is trivial on $N$, equivalently, that
the restriction $M=L|_Y$ is trivial on $Y$. 
We can assume by Remark \ref{(2)} that $L=\pi^*(M)$ for the 
projection $\pi : N\to Y$. 
Let $N^o=N-Y$ and $M^o=M-Y$, and let
$\nu : N^o\to Y$ and $\mu : M^o\to Y$ denote the projections. 
We have the commutative diagram of homotopy exact sequences
$$\begin{array}{cccccccc}
\pi_2(Y) & \to &\pi_1(N^o_y) & \to & \pi_1(N^o) & \stackrel{\nu_*}{\to} & 
\pi_1(Y) & \to 1\\
 \|&& \dna & &             \dna & &            \| & \\
\pi_2(Y) & \to &\pi_1(M^o_y) & \to & \pi_1(M^o) & \stackrel{\mu_*}{\to} & 
\pi_1(Y) & \to 1
\end{array}$$
where the middle two vertical maps are induced by $\sigma$.
As $|\nu_Y|(\sigma)=0$, the map $\pi_1(N^o_y)\to \pi_1(M^o_y)$
is zero. Hence from the above diagram there exists a group
homomorphism $\psi : \pi_1(Y) \to \pi_1(M^o)$ such that 
$\mu_*\circ \psi = \id_{\pi_1(Y)}$. 

As complex line bundles $M$ are classified by their first Chern
class $c_1(M)\in H^2(Y;\Z)$, we have to show that $c_1(M)=0$.
Given any non-zero element $c$ of $H^2(Y;\Z)$, there exists a compact
$2$-manifold $K$ and a continuous map $f : K \to Y$, such that
$f^*(c) \ne 0$ in $H^2(K;\Z)$. The section 
$f^*(\sigma) \in \Gamma(f^*N,f^*L)$ has the same vanishing multiplicity
as $\sigma$. Hence we can base change to $K$, and therefore
we can assume that $Y$ is a compact
$2$-manifold for the sake of proving that $L$ is trivial.

Unless such a $Y$ is homeomorphic to $S^2$ or $\PP_{\R}^2$
(we treat these two cases separately later), 
we know that $Y$ is a $K(\pi, 1)$ space. 
Hence it follows from the homotopy exact
sequence for $\mu : M^o\to Y$ that 
$M^o$ is also a $K(\pi,1)$ space. Hence the group
homomorphism $\psi : \pi_1(Y) \to \pi_1(M^o)$ is induced
by a map $s : Y \to M^o$. As $\mu_*\circ s_* = \id_{\pi_1(Y)}$,
and as $Y$ is a $K(\pi,1)$ space, the map $\mu\circ s : Y\to Y$
is homotopic to $\id_Y$. Hence by the homotopy lifting
property of $\mu : M^o\to Y$, the map $\id_Y : Y\to Y$
lifts to give a section of $\mu$, proving triviality of $M$.

Next, suppose $Y = S^2$. Then $M= \O(d)$ for  
$d \in Z$, with $c_1(M) = d \in H^2(Y;\Z)=\Z$. It can be seen that
$\pi_1(M^o) = \Z/(d)$. Now 
as the map $\pi_1(N^o_y)\to \pi_1(M^o_y)$
is zero, from the above commutative diagram
comparing the long exact homotopy sequences we 
see that the map $\pi_2(Y) \to \pi_1(M^o_y)$
is zero, hence the map
$\pi_1(M^o_y) \to \pi_1(M^o)$ is an isomorphism, showing
$\pi_1(M^o) = \Z$, therefore $M$ is trivial when $Y = S^2$.

Finally, suppose $Y = \PP_{\R}^2$. Then as $H^2(Y;\Z) = \Z/(2)$,
$Y$ has exactly one non-trivial complex line bundle $M$, up to isomorphism.
It can be seen that for this $M$, we have $\pi_1(M^o) = \Z$,
so there cannot exist a section to $\pi_1(M^o)\to \pi_1(Y) = \Z/(2)$.
Hence $M$ must be trivial, even when $Y = \PP_{\R}^2$. 

Now choose a non-vanishing global section $e$ of $L$,
and let $g : \pa D \to \C\,^*$ be defined by 
$\sigma = g e$ on $\pa D$. The hypothesis $\nu_Y(\sigma)=0$
implies that the restriction of $g$ to any fiber of $p: \pa D \to Y$
is homotopic to the constant function $1$. 
By Lemma \ref{BZ} above, 
$g$ admits a prolongation $\ov{g} : D \to \C\,^*$.
Put $\sigma' = \ov{g}e$ on $D$, and $\sigma' = \sigma$ on $N - D$. 
This proves the lemma. \hfill$\square$.

\begin{lemma}\label{L is power of normal} 
Let $Y$ be a connected Hausdorff space of 
homotopy type of a CW complex, and
let $\pi : N\to Y$ be a real vector bundle on $Y$ of rank $2$,
with a metric. 
Let $M \to Y$ be a complex line bundle on $Y$, let 
$L = \pi^*(M)$, and let
$u \in \Gamma(N,L)$ such that $u$ is non-vanishing on $N-Y$
and $\nu_Y(u)=m\ge 1$. Let $N$ be given the $\C$-linear structure induced
by the orientation defined by $u$ (see Remark \ref{(4)}), 
making it a complex
line bundle on $Y$. Let $N^m$ be the $m$ th tensor power of
the complex line bundle $N$.

Then we have the following.

(i) The complex line bundle $M$ is isomorphic to
$N^m$. 

(ii) In fact, there exists an isomorphism $\psi : N^m \to M$ of
complex line bundles on $Y$ with the following property. Let 
$\tau \in \Gamma(N,\pi^*(N))$ be the tautological section
(defined by $\id : N\to N$), and let $\tau^m\in 
\Gamma(N,\pi^*(N^m))$ be its $m$ th power.
Consider the section $\psi(\tau^m) \in \Gamma(N,L)$, 
and let $f: N-Y \to \C\,^*$
be defined by $u|_{N-Y} = f\cdot \psi(\tau^m)|_{N-Y}$. 
Then $f$ is homotopic to the constant function $1$. 
\end{lemma}

\proof As $\nu_Y(u)\ne 0$, by Remark \ref{(4)} $u$ defines 
an orientation on the real vector bundle $N$. Together with the 
chosen  metric on $N$, this gives a $\C$-linear structure on 
$N$, so we regard it as a complex line bundle. Let 
$K$ be the complex line bundle on $Y$ defined by
$K = M\otimes N^{-m}$.
Let $\pi^*(K)$ denote the pullback of $K$ to the total space of $N$.
Let the continuous map $v : N \to M$ 
be the composite of the section $u: N \to N\times_YM$ with
the projection $N\times_YM\to M$. 
(Note that $v$ is over $Y$, but may not be linear on fibers.)
Given any $x\in N-Y$ over $y\in Y$, let 
$$\sigma_x : (N_y)^{\otimes m} \to M_y$$
be the unique $\C$-linear map under which $x^{\otimes m}\mapsto v(x)$.
We can regard $\sigma_x$ as an element of $K_y$, and as
$x$ varies we get a section $\sigma \in \Gamma(N-Y,\pi^*(K))$.
From the hypothesis that $\nu_Y(u)=m$ it follows that 
$\nu_Y(\sigma)=0$. Hence by Lemma 
\ref{non-vanishing prolongation}, the complex line
bundle $K$ is trivial. This proves (i).

Now suppose $\psi_1 : N\to M$ is some isomorphism of complex bundles,
which exists by (i).
Consider the resulting section 
$\psi_1(\tau^m) \in \Gamma(N,M_N)$, and let $f_1: N-Y \to \C\,^*$
be defined by $u = f_1\psi_1(\tau^m)$. As both $u$ and $\psi_1(\tau^m)$ have 
oriented topological vanishing multiplicity $1$ along $Y$,
it follows that $f_1$ has topological vanishing multiplicity $0$ along $Y$.
By Lemma \ref{BZ}, there exists a function
$f_2 : N \to \C\,^*$ such that $f_2|_{N-D} = f_1|_{N-D}$ 
where $D\subset N$ is the unit disk bundle in $N$. 
Now put $\psi = f_2\psi_1(\tau^m)$. Then 
$u|_{N-Y} = f \cdot \psi|_{N-Y}$, where 
$f = f_1/f_2 : N-Y \to \C\,^*$. 
As $f|_{N-D}$ is the constant function $1$, and as
$N-D\hra N-Y$ is a homotopy equivalence, $f|_{N-Y}$ is homotopic
to the constant function $1$. This proves (ii) and completes
the proof of the lemma. \hfill $\square$

%
%

\section{Proof of the Main Theorem}

In this section we allow the singular cohomologies to have
coefficients in an arbitrary ring $R$.

\centerline{\bf Vanishing multiplicity of discriminant}

Recall that for any triple $(E,L,b)$ on $X$, 
the form $b$ is given on a trivializing open cover of $X$ 
by $r\times r$ matrices of functions (where $r=\rank(E)$), 
whose determinants fit together to define a global section
$\det(b) \in \Gamma(X,  L^r\otimes \det(E)^{-2})$
called the {\bf discriminant} of $(E,L,b)$.
It is clear that if $K$ is any line bundle on $X$, then for the triple
$(E\otimes K, L\otimes K^{\otimes 2}, b\otimes 1_K)$ we have  
$\det(b\otimes 1_K) = \det(b) \in  \Gamma(X,  L^r\otimes \det(E)^{-2})$.  
Hence for a quadric bundle $Q=[E,L,b]$, the section
$\det(Q) = \det(b)$ is well defined.

The non-degeneracy condition on $(E,L,b)$ just means that
$\det(b)$ is nowhere zero. More generally, 
if $Z\subset X$ is the vanishing locus of $\det(b)$ then the 
restriction $(E,L,b)|_{X-Z}$ is non-degenerate. 

Now consider a pair $(X,Y)$ where $X$ is a topological space, and $Y$ a
topological divisor in $X$. Given a quadric bundle  $Q =[E,L,b]$ 
on $X$ which is
non-degenerate on $X-Y$, the non-negative integer valued
functions $|\nu_Y|(\det(Q))$ and $\ovnu_Y(\det(Q))$ on $Y$ are the 
absolute vanishing
multiplicity of the discriminant $\det(b)$ along $Y$, and  
its parity, as in 
Definition \ref{(3)} and Remark \ref{tubular neighbourhood} above.

\medskip

\medskip

\centerline{\bf Minimally degenerate quadric bundles and canonical triples}

Recall from the introduction that a triple $\T = (E,L,b)$,
or the corresponding quadric bundle $Q=[E,L,b]$ on a base $Y$ is 
{\bf minimally degenerate} if $\rank(b_y) = \rank(E) -1$
at all points $y\in Y$. 
This means that the linear map $b : E \to L\otimes E^*$
is of constant rank, with kernel a line subbundle $K\subset E$. 
Let $\ov{E} = E/K$ be the quotient bundle, on which we get an induced
bilinear form $\ov{b} : \ov{E}\otimes \ov{E} \to L$
which is non-degenerate. Hence we get a non-degenerate triple
$\ov{\T} = (\ov{E},L,\ov{b})$
of rank $r-1$ where $r=\rank(E)$.
Recall that we have defined the {\bf canonical triple} 
$\T^Q$ associated to such a minimally degenerate quadric bundle $Q$ 
to be the non-degenerate rank $r-1$ triple
$$\T^Q = \ov{\T} \otimes \ker(b)^{-1} =
(\ov{E}\otimes K^{-1}, L\otimes K^{-2}, \ov{b} \otimes 1_{K^{-1}})$$
This is well-defined, independent of the 
choice of the representative $\T$ for
$Q$. It represents the non-degenerate quadric bundle 
$\ov{Q} = [\ov{\T}]$ of rank $r-1$.

\medskip

\centerline{\bf Mild degeneration}

Let $Y\subset X$ be a topological divisor.
Recall from the introduction that we call $\T=(E,L,b)$ 
a {\bf mildly degenerating triple} on $(X,Y)$ if 
$\T_{X-Y}$ is non-degenerate of rank $r$ say, and
$\T_Y$ is of constant rank $r-1$. We say that the corresponding
$Q = [\T]$ is a {\bf mildly degenerating quadric bundle} on $(X,Y)$
of {\bf generic rank} $r$.

In the algebraic category, we have the following.

\begin{proposition}\label{Q nonsingular} 
Let $X$ be a non-singular complex algebraic variety, 
and let $Y\subset X$ be
a smooth divisor. Let $Q\to X$ be mildly degenerating 
quadric bundle over $(X,Y)$ in the algebraic category. 
Then the following are equivalent.

(i) The total space of $Q$ is non-singular.

(ii) The algebraic vanishing multiplicity $\nu_Y(\det(Q))$ 
(which is the same as the oriented topological vanishing multiplicity 
with respect to the natural orientation) is $1$ over 
all components of $Y$.
\end{proposition}

\proof This is a simple generalization of the Proposition 3 of [N],
and follows by the Jacobian criterion.
\hfill$\square$

\medskip

\centerline{\bf Orthogonal decomposition} 

\defin\label{triple sum} 
{\bf (Direct sums of triples) } 
Let $(E_1,L,b_1)$ and $(E_2,L,b_2)$ be two triples on a space $X$,
where the line bundle $L$ is common to the two triples. 
We define the {\bf direct sum} of $(E_1,L,b_1)\oplus (E_2,L,b_2)$ 
to be the triple $(E_1\oplus E_2,\,L,\,b_1\oplus b_2)$, where
$b_1\oplus b_2$ has the matrix form 
$\left(\begin{array}{cc}  
b_1 & \\ 
    & b_2 
\end{array}\right)$. 
More generally, if $(E_1,L_1,b_1)$ and $(E_2,L_2,b_2)$ are
two triples, and we are given an isomorphism
$\psi : L_1\to L_2$, Then the triple 
$(E_1,L,b_1)\oplus_{\psi} (E_2,L,b_2)$, which we call the
{\bf direct sum via $\psi$}, is defined to be
$(E_1\oplus E_2,\,L_2,\,\psi b_1\oplus b_2)$, where
$\psi b_1\oplus b_2$ has the matrix form 
$\left(\begin{array}{cc}  
\psi b_1 & \\ 
         & b_2 
\end{array}\right)$.

\begin{lemma}\label{orthogonal decomposition} 
Let $\pi: N\to Y$ be a rank $2$ real vector bundle, 
and let $(E,L,b)$ be a mildly degenerating triple on $(N,Y)$,
of generic rank $r$. Let the line subbundle $K \subset E_Y$ be the 
kernel of the restriction of $b$ to $Y$, and 
let $(E_Y/K, L_Y, \ov{b_Y})$ be the non-degenerate
triple of rank $r-1$ on $Y$, where $\ov{b_Y}$ is induced by $b$.
Then there exists a bilinear form 
$b' : \pi^*(K)\otimes \pi^*(K) \to \pi^*(L_Y)$ such that
the triple $(E,L,b)$ on $N$ is isomorphic to the direct sum
$$(E,L,b)\cong (\pi^*(K), \pi^*(L_Y),b') \oplus 
\pi^*(E_Y/K, L_Y, \ov{b_Y})$$
\end{lemma}

\proof As $Y\subset N$ is a strong deformation retract,
we assume without loss of generality that $E$ is of the form
$\pi^*(E_Y)$ and $L$ is of the form $\pi^*(L_Y)$.
We choose a Hermitian metric $h_Y$ on the vector bundle $E_Y$.
Let $E''_Y$ be the orthogonal complement of $K$ in $E_Y$
with respect to $h_Y$.
Then the quotient map $E_Y \to E_Y/K$ gives an 
isomorphism $E''_Y\to E_Y/K$. 

Let $E'' = \pi^*(E''_Y)$ which is a rank $n-1$ 
subbundle of $E = \pi^*(E_Y)$. 
Let $b'' : E''\otimes E'' \to L$ be induced by 
$b : E\otimes E \to L$. Note that $b''|_Y = \ov{b_Y}$, so
the discriminant $\det(b'')$ is non-zero on $Y$. 
Since $Y$ is paracompact we can find a positive continuous function $f$ on 
$Y$ such that for any $x\in N$ we have
$\det(b'')_x \ne 0$ whenever $0\le |\!|x|\!| \le f(\pi (x))$.
By scaling the metric using $f$, we assume that 
$\det(b'')_x \ne 0$ whenever $0\le |\!|x|\!| \le 1$.

Consider the continuous map $\gamma : N \to N$ defined by
$$\gamma(x) = \left\{ \begin{array}{ll}
x & \mbox{ if } |\!|x|\!| \le 1,  \\
x / |\!|x|\!| & \mbox{ if } |\!|x|\!| \ge 1.  
\end{array}\right.$$
The triple $(E,L,b)$ on $N$ is then isomorphic to
the pullback $\gamma^*(E,L,b)$, which allows us to assume 
without loss of generality that
the restriction $b''= b|_{E''}$ is everywhere non-degenerate on $N$.

Let $E'\subset E$
be the orthogonal complement of $E''$ with respect to
the form $b : E\otimes E \to L$. 
As $b$ is non-degenerate on $E''$, it follows that
$E'$ is a line subbundle of $E$. 
From its definition, we have the following
equality of subbundles of $E_Y$
$$E'|_Y = K \subset E_Y$$
Let $b'=b|_{E'}$. By its construction, the triple
$(E',L,b')$ on $N$ is non-degenerate of rank $1$ on $N-Y$,
and degenerates on $Y$, and we have the direct sum
decomposition
$$(E,L,b) = (E',L,b') \oplus (E'',L,b'')$$
Note that as $(E'',L,b'')$ is non-degenerate, it is equivalent to
a principal $GO(n-1)$-bundle $P$ on $N$. As 
$Y\subset N$ is a strong 
deformation retract with inverse $\pi : N\to Y$, 
it follows that $P$ is isomorphic to $\pi^*(P|_Y)$. 
Hence $(E'',L,b'')$ is isomorphic to $\pi^*(E_Y/K, L_Y, \ov{b_Y})$.
Also note that for the same reason, 
the line bundle $E'$ is isomorphic to $\pi^*K$.
This completes the proof of the lemma. \hfill$\square$

\medskip

\centerline{\bf Non-degenerate prolongation when $|\nu_Y|(\det(b))=0$}

\begin{lemma}\label{non-degenerate prolongation} 
If the absolute vanishing multiplicity 
$|\nu_Y|(\det(b))=0$, then there exists an everywhere 
non-degenerate symmetric bilinear form $b_* : E\otimes E \to L$
such that $b_*$ coincides with $b$ on 
$N^{\ge 1}  = \{ v\in N \,|\, |\!| v|\!| \ge 1\}$. 
In particular, all quadric invariants $\alpha \in H^*(X-Y)$
of $[E,L,b]_{X-Y}$ prolong to $H^*(X)$ as the corresponding
invariants of $[E,L,b_*]$, hence map to $0$ under the 
Gysin boundary map $H^*(X-Y)\to H^{*-1}(Y)$.
\end{lemma} 

\proof By Lemma \ref{orthogonal decomposition}, we can assume that
$$(E,L,b) = (\pi^*(K), \pi^*(L_Y),b') \oplus 
\pi^*(E_Y/K, L_Y, \ov{b_Y})$$
As $\ov{b_Y}$ is non-degenerate, it follows that 
$|\nu_Y|(\det(b'))= |\nu_Y|(\det(b))=0$.
We now apply Lemma \ref{non-vanishing prolongation} with
$L = \pi^*(L_Y\otimes K^{-2})$ and $\sigma =\det(b')$.
Hence there exists a non-vanishing
section $b'_* \in \Gamma(N, \pi^*(L_Y\otimes K^{-2}))$ such that
$b'_*$ coincides with $b'$ on $N-D$, where $D\subset N$
is the unit disk subbundle. Consider the resulting 
triple $(\pi^*(K), \pi^*(L_Y),b'_*)$. Then the direct sum
$(\pi^*(K), \pi^*(L_Y),b'_*)\oplus \pi^*(E_Y/K, L_Y, \ov{b_Y})$ 
satisfies the lemma.
\hfill$\square$

\medskip

\centerline{\bf Standard models of mildly degenerating triples}

In view of the above lemma, 
to prove the main theorem we only need to consider
the situation where the absolute vanishing multiplicity 
$|\nu_Y|(\det(Q)) \ge 1$. Therefore, $N$ has a unique
orientation such that the oriented vanishing multiplicity 
$\nu_Y(\det(Q)) \ge 1$. Together with the metric
on $N$, this defines a $\C$-linear structure on $N$.
Hence in what follows, we assume that $N$ is a complex line bundle, and
$\nu_Y(\det(Q))= m \ge 1$.

Now let $\tau \in \Gamma(N, \pi^*(N))$
be the tautological section defined by the identity 
homomorphism $N\to N$. Note that its $m$ th tensor power
$\tau^m \in  \Gamma(N, \pi^*(N^m))$, where
$N^m$ is the $m$ th tensor power of the line bundle $N$,
is non-vanishing outside $Y$ and has oriented topological vanishing
multiplicity $m$ along $Y$. On the total space of $N$, consider the triple
$(\O_N,\pi^*(N^m), \tau^m)$
where the bilinear form $\tau^m : \O_N\otimes \O_N \to\pi^*(N^m)$ 
sends $1\otimes 1 \mapsto\tau^m \in  \Gamma(N, \pi^*(N^m))$.
Note that this triple is non-degenerate of rank $1$ outside 
$Y$ and minimally degenerate on $Y$, with $\nu_Y(\det(\tau^m))=m$.

Now let $(F,N^m,q)$ be a non-degenerate triple of rank  
$r-1$ on $Y$, where $r\ge 2$. 
Then on $N$ we have the direct sum triple
$$\M = (\O_N,\pi^*(N^m), \tau^m)\oplus (F,N^m,q)$$
From its construction, the above triple  
is mildly degenerating on $(N,Y)$, with $\nu_Y(\det(\M)) = m$.

We call this triple $\M$ as
the {\bf standard model} of a mildly degenerating triple on $(N,Y)$,
corresponding to the data $(F,N^m,q),\,m$.
This name is justified by the Remark \ref{model vindicated}  
which follows the lemma below.

\begin{lemma}\label{standard model} 
Let $N\to Y$ be a complex line bundle, 
and let $(E,L,b)$ be a mildly degenerating triple on $(N,Y)$
of generic rank $r$, with oriented vanishing multiplicity 
$\nu_Y(\det(b))= m \ge 1$. 
Let $K \subset E_Y$ be the kernel of 
$b_Y : E_Y \to L_Y \otimes E_Y^*$, and 
let $(E_Y/K, L_Y, \ov{b_Y})$ be the induced non-degenerate
triple of rank $r-1$ on $Y$. 
Then there exists an isomorphism $\psi: N^m \otimes K^2\to L_Y$ 
of complex line bundles, such that
on $N-Y$,  we have an isomorphism of non-degenerate triples
$$(E,L,b)_{N-Y} \cong (\O_{N-Y},\pi^*(N^m), \tau^m)\otimes \pi^*(K)\,
\oplus_{\psi} \, \pi^*(E_Y/K, L_Y, \ov{b_Y})$$
where $N^m\otimes K^2$ is identified with $L_Y$ via $\psi$ for
defining the direct sum of triples on the right hand side as
in Definition \ref{triple sum}, and $\pi: N-Y \to Y$ is the projection. 
\end{lemma} 

\proof By Lemma \ref{orthogonal decomposition}, we can assume that 
$$(E,L,b) = (\pi^*(K), \pi^*(L_Y),b') \oplus 
\pi^*(E_Y/K, L_Y, \ov{b_Y})$$ 
As $\ov{b_Y}$ is non-degenerate, 
it follows that $\nu_Y(\det(b'))= \nu_Y(\det(b))=m \ge 1$.  
We now apply Lemma \ref{L is power of normal} with 
$M = L_Y\otimes K^{-2}$ and $u=\det(b')$.
Let $\psi: N^m \to L_Y\otimes K^{-2}$ 
be an isomorphism of complex line bundles as 
given by Lemma \ref{L is power of normal}(ii), so 
that 
$$\det(b')|_{N-Y} = f \cdot \psi(\tau^m)|_{N-Y}$$
where $f : N-Y \to \C\,^*$ is homotopic to the constant map
$1: N-Y \to \C\,^*$. Hence $f$ admits a continuous square-root
$g : N-Y \to \C\,^*$, with $g^2=f$.
It follows that the scalar multiplication
$g : K\to K$ gives an 
isomorphism of triples
$$(g,\id_L) : (\pi^*(K), \pi^*(L_Y), b') \to
(\pi^*(K), \pi^*(L_Y), \psi(\tau^m))$$
Substituting this in the direct sum decomposition
of $(E,L,b)$ at the beginning of the proof, the lemma follows.
\hfill$\square$ 

\rem\label{model vindicated}  
With notation as in the lemma above,
let $F$ be the vector bundle $(E_Y/K)\otimes K^{-1}$ on $Y$. 
Then under the isomorphism $\psi : N^m \otimes K^2 \to L_Y$,
the triple $(E_Y/K, L_Y,  \ov{b_Y})$ takes the form
$(F,N^m,q)\otimes K$. Hence the above lemma says that 
given any triple $\T =(E,L,b)$ which minimally degenerates on $Y$,
there exists a standard model triple 
$\M = (\O_N,\pi^*(N^m), \tau^m)\oplus \pi^*(F,N^m,q)$
with the following properties. 

(1) The triples $\M\otimes K$ and $\T$ have isomorphic restriction
on $Y$. 

(2) The triples $\M\otimes K$ and $\T$ have isomorphic restriction
on $N-Y$. 

(3) The oriented topological vanishing multiplicity 
$\nu_Y(\det(\M)) = \nu_Y(\det(\M\otimes K))$ of $\M$
equals the oriented topological vanishing multiplicity 
$\nu_Y(\det(\T))$ of $\T$.

Hence for the purpose of studying the topology of degeneration
where $\nu_Y(\det(b)) \ge 1$,
we can confine ourselves to quadric bundles $[\M]$ defined by
standard models $\M$.

\bigskip
\pagebreak

\centerline{\bf Proof of the Main Theorem for $|\nu|_Y(\det(Q))$ even}


Let $N\to Y$ be a complex line bundle. Consider a standard model triple 
$$\M = (\O_N,\pi^*(N^m), \tau^m)\oplus \pi^*(F,N^m,q)$$ 
on $N$, where $m = 2k>0$ is even. We define a new triple $\N$ by
$$\N  = (\pi^*(N^k),\pi^*(N^m), t)\oplus \pi^*(F,N^m,q)$$
where $t : \pi^*(N^k)\otimes \pi^*(N^k) \to \pi^*(N^{2k})
=\pi^*(N^{m})$ is the tensor product.
We have a morphism of triples
$f = (\tau^k,\id) : (\O_N,\pi^*(N^m), \tau^m) \to 
(\pi^*(N^k),\pi^*(N^m), t)$, which together with 
identity on $\pi^*(F,N^m,q)$ gives a morphism $g : \M \to \N$.
From its definition, it is clear that $g|_{N-Y}$ is an isomorphism
of triples on $N-Y$, and also $g|_Y$ is an isomorphism of
triples on $Y$.

We have thus shown that for any mildly degenerating triple $\T$ 
on $(N,Y)$ where $N$ is a rank $2$ real vector bundle on $Y$, 
such that $|\nu|_Y(\det(\T))$ is even, 
there exists an everywhere non-degenerate triple 
$\T'$ on $N$ such that $\T$ is isomorphic to
$\T'$ on $N-Y$.
In particular, all quadric invariants 
$\alpha[\T_{N-Y}] \in H^*(N-Y)$
of $[\T]_{N-Y}$ prolong to $H^*(N)$ as the corresponding
invariants of $[\T']$, hence map to $0$ under the 
Gysin boundary map $H^*(N-Y)\to H^{*-1}(Y)$.
This completes the proof of the main theorem when the topological
vanishing parity $\ovnu_Y(\det(Q))$ is zero. \hfill$\square$

\medskip

\centerline{\bf Proof of the Main Theorem for $|\nu|_Y(\det(Q))$ odd} 

Let $N\to Y$ be a complex line bundle. Consider a standard model triple 
$$\M = (\O_N,\pi^*(N^m), \tau^m)\oplus \pi^*(F,N^m,q)$$ 
on $N$, where $m = 2k+1>0$ is odd. We define a new triple $\N$ by
$$\N  = (\pi^*(N^k),\pi^*(N^m), \tau\circ t)\oplus \pi^*(F,N^m,q)$$
where $t : \pi^*(N^k)\otimes \pi^*(N^k) \to \pi^*(N^{2k})$ is the 
tensor product, and $\tau : \pi^*(N^{2k})\to \pi^*(N^{2k+1})
=\pi^*(N^{m})$ is induced by $\tau \in \Gamma(N,\pi^*(N))$.  
Again, we have a morphism of triples
$f = (\tau^k,\id) : (\O_N,\pi^*(N^m), \tau^m) \to 
(\pi^*(N^k),\pi^*(N^m), t)$, which together with 
identity on $\pi^*(F,N^m,q)$ gives a morphism
$g: \M \to \N$. 
such that $g|_{N-Y}$ is an isomorphism
of triples on $N-Y$, and also $g|_Y$ is an isomorphism
of triples on $Y$.
Moreover, note that $\nu_Y(\det(b_{\N}))=1$.
The triple $\N \otimes \pi^*(N^{-k})$ is the standard model triple
$$\MM = (\O_N, \pi^*(N), \tau) \oplus 
\pi^*(F\otimes N^{-k}, N, q\otimes 1_{N^{-k}})$$
with $\nu_Y(\det(\MM))=1$. Hence we have the following.

\begin{lemma}\label{odd=1} 
Let $\pi : N\to Y$ be a complex line
bundle. Given any standard model triple $\M$ on $N$ 
with oriented topological vanishing multiplicity 
$\nu_Y(\det(\M)) = 2k + 1>0$, 
there exists a standard model triple $\MM$ on $N$ 
with the following three properties:  
(i) the triples $\M\otimes \pi^*(N^{-k})$ and $\MM$ 
have isomorphic restriction on $Y$,  
(ii) the triples $\M\otimes \pi^*(N^{-k})$ 
and $\MM$ have isomorphic restriction
on $N-Y$, and 
(iii) The oriented topological vanishing  multiplicity of $\MM$ is $1$.
\end{lemma}

Hence to prove the main theorem,
it is enough to restrict ourselves to the case of
quadric bundles $[\MM]$ on $N$, defined by triples of the form 
$$\MM = (\O_N, \pi^*(N), \tau) \oplus \pi^*(V,N,b)$$
where $(V,N,b)$ is a non-degenerate triple on $Y$.
Note that for such an $[\MM]$, the induced canonical triple
$\T^{[\MM]}$ on $Y$ is just $(V,N,b)$.

\medskip

\centerline{\bf Pairs of classifying maps}

Let $\T =(V,N,b)$ be a rank $m$ non-degenerate triple on $Y$. 
On $N_o= N-Y$, the restriction of the tautological section
$\tau : \O_N \to \pi^*(N)$ admits an inverse $t_{N_o}$. 
We get a non-degenerate pair 
$t_{N_o}(\T)= (\pi^*(V),t_{N_o}\circ \pi^*(b))$
of rank $m$ on $N-Y$ (which means $t_{N_o}\circ \pi^*(b)$
is an $\O_{N-Y}$-valued non-degenerate symmetric bilinear form
of rank $m$ on $\pi^*(V)$).
Let $BGO(m)$ be a classifying space of $GO(m)$, where $m\ge 1$, and 
let $\U = (\EE,\LL,\bb)$ denote the universal triple on $BGO(m)$.
Consider the non-degenerate pair 
$t_{\LL}(\U)$ on $\LL_o = \LL - BGO(m)$, 
which gives an identification of $\LL_o$ with 
the classifying space $BO(m)$ of $O(m)$.
Let $\eta: Y\to BGO(m)$ be a classifying map
for $(V,N,b)$. 
By the defining property of a classifying map, there exists  
an isomorphism of
triples $(\alpha,\beta): \T \to \eta^*\U$
where $\alpha : V \to \eta^*\EE$ and
$\beta : N \to \eta^*\LL$ are vector bundle isomorphisms which take
$b$ to $\eta^*\bb$. Let $\eta': \eta^*\LL \to \LL$ be the projection map,
and consider the composite map $\eta'\circ \beta : N\to \LL $.
Then from its construction, its restriction 
$\theta =(\eta'\circ \beta)|_{N_o} :  N_o \to \LL_o = BO(m)$
is a classifying map for the non-degenerate 
quadratic pair 
$t_{N_o}(\T) = (\pi^*(V),\, t_{N_o}\circ \pi^*(b))$ on $N_o$.
Moreover, we have a pullback diagram (Cartesian square) of vector
bundles on $Y$ as follows 
$$\begin{array}{ccc}
N & \stackrel{\eta'\circ \beta}{\lra} & \LL \\
\dna &      & \dna \\
Y & \stackrel{\eta}{\lra} & BGO(m)
\end{array}$$

As the Gysin boundary maps $\pa$ commute 
with pullback of bundles, the above proves
the following crucial lemma, by taking 
$\MM = (\O_N,\pi^*(N),\tau)\oplus \pi^*(V,N,b)$,
and choosing $\theta$ as above.

\begin{lemma}\label{model Gysin}  
Let $\pi: N\to Y$ be a complex line bundle,
let $\T = (V, N,b)$ be a non-degenerate triple of rank 
$n-1$ on $Y$, where $n\ge 2$, and let $\MM$ be the mildly degenerating
triple
$$\MM = (\O_N, \pi^*(N), \tau) \oplus \pi^*(V,N,b)$$
on $(N,Y)$. Let $\eta : Y\to BGO(n-1)$ be a classifying map for
the non-degenerate triple $\T^{[\MM]_Y} = \T$, and  
let $\theta : N-Y \to BO(n-1)$ be a classifying map for
the principal $O(n-1)$-bundle $t_{N_o}(\T)$. 
Then the following diagram is commutative
for any integer $i$.
$$\begin{array}{ccc}
H^i(BO(n-1)) & \stackrel{\theta^*}{\to} & H^i(N-Y) \\ 
{\scriptstyle\pa} \dna ~~&      &~~ \dna {\scriptstyle\pa}\\
H^{i-1}(BGO(n-1))  & \stackrel{\eta^*}{\to} & H^{i-1}(Y)
\end{array}$$
\end{lemma}

The above lemma shows that 
Gysin boundary acts on invariants $\alpha[\MM_{N-Y}]$
of a mildly degenerating quadric $\MM$ on $(N,Y)$
defined by a triple $(\O_N, \pi^*(N), \tau) \oplus \pi^*(V,N,b)$
exactly as given by the main theorem. As we have already 
reduced the proof of the main theorem to triples of the above form, 
this completes the proof of the main theorem.
\hfill$\square$

%
%

\section{Sample calculations}

The main theorem, combined with our earlier determination of
the Gysin boundary for the fibrations $BO(2n+1)\to BGO(2n+1)$ and
$BO(2n) \to BGO(2n)$
gives the following corollaries which tell us how to explicitly calculate
the Gysin boundary image of any given quadric invariant.

\begin{corollary}\label{odd to even} 
{\rm (Odd rank degenerating to even rank) }
Let $Q$ be a mildly degenerating quadric bundle on $(X,Y)$,
of generic rank $2n+1$, where $Y$ is connected. 
Suppose that the degeneration multiplicity $|\nu|_Y(\det(Q))$ is odd. 
Then under the Gysin boundary map
$\pa : H^*(X-Y)\to H^{*-1}(Y)$, 
quadric invariants $h(\wh{w}_2,\,\ldots,\,\wh{w}_{2n+1})$ of $Q_{X-Y}$ 
map to characteristic classes  $g(\lambda,\,a_{2i-1},\,b_{4j},\,d_T)$
of the triple $\T^{Q_Y}$ on $Y$ as follows.
Under the substitutions 
$\wh{w}_r \mapsto 
{2n \choose r}\,w_1^r + \sum_{2\le i \le r} 
{2n+1-i \choose r-i}\, w_1^{r-i}\,w_i$
for $2\le r \le 2n+1$ followed by putting $w_{2n+1}=0$,  
let the polynomial $h(\wh{w}_2,\,\ldots,\,\wh{w}_{2n+1})$
map to the sum $\sum_T P_T$ indexed by subsets 
$T\subset \{ 1,\ldots,n\}$, where for each 
$T = \{ i_1 <\ldots < i_r \}$, $P_T$ has the form
$$P_T = f_T(w_{2i-1})\cdot g_T(w_{2i}^2)\cdot 
w_{2i_1} \cdots w_{2i_r}$$
where $f_T$ and $g_T$ are polynomials in $n$ variables each.
Then $\pa (h) = \sum_T \pa (P_T)$ where 
$$\pa (P_T) = 
\left\{ 
\begin{array}{ll}
0 & \mbox{if } r= 0 \\
f(a_{2i-1})\cdot g(b_{4i})
\cdot\, a_{2i_1-1} & \mbox{if } r=1 \\
f(a_{2i-1})\cdot g(b_{4i})
\cdot d_T & \mbox{if } r \ge 2, ~\mbox{where } \\
& T = \{ i_1 <\ldots < i_r \}.
\end{array}\right.
$$
If the degeneration multiplicity $|\nu|_Y(\det(Q))$ is 
even, then $\pa(h)=0$.
\end{corollary}

\begin{corollary}\label{even to odd} 
{\rm (Even rank degenerating to odd rank) }
Let $Q$ be a mildly degenerating quadric bundle on $(X,Y)$ 
of generic rank $2n+2$, where $Y$ is connected. 
If the degeneration multiplicity $|\nu|_Y(\det(Q))$ is odd,
then under the Gysin boundary map
$\pa : H^*(X-Y)\to H^{*-1}(Y)$, a quadric invariant of the form
$h(\lambda,\,a_{2i-1},\,b_{4j},\,d_T)$ of $Q_{X-Y}$
maps to the characteristic class 
$g(c,\wh{w}_2,\,\ldots,\,\wh{w}_{2n+1})$ of
the triple $\T^{Q_Y}$ as follows. 
Under the substitution 
$$\lambda\mapsto 0,~
a_{2i-1}\mapsto w_{2i-1},~  
b_{4j}\mapsto w_{2j}^2 \mbox{ and } 
d_T\mapsto \sum_{i\in T} w_{2i-1}v_{T-\{ i\}}$$ 
where $w_{2n+2}=0$, let the polynomial 
$h(\lambda,\,a_{2i-1},\,b_{4j},\,d_T)$ map to 
the polynomial 
$h'(w_1,\ldots,w_{2n+1})$. 
Under the change of variables from 
$ w_1,\ldots,w_{2n+1}$ to $w,\,\wh{w}_2,\,\ldots$,
$\wh{w}_{2n+1}$ given by
Lemma \ref{w in terms of wbar}, let $h'$ map to 
the sum $\sum_iw^if_i(\wh{w}_2,\ldots,\wh{w}_{2n+1})$ 
where $f_i$ are polynomials in $2n$ variables.
We have 
$$\pa(h) =  
\sum_{j\ge 0} c^j f_{2j+1}(\wh{w}_2,\ldots,\wh{w}_{2n+1})$$ 
If the degeneration multiplicity $|\nu|_Y(\det(Q))$ is 
even, then $\pa(h) =0$.
\end{corollary}

\rem {\bf Even to odd degeneration over $(N,Y)$ } 
Let $\pi : N\to Y$ be a rank $2$ real vector
bundle, and $\T = (E,L,b)$ a mildly degenerating triple 
on $(N,Y)$ of generic rank $2n$. 
Consider the resulting principal $GO(2n-1)$-bundle
$\ov{\T_Y}$ on $Y$. 
The character $\kappa : GO(2n-1) \to \C^*$ defined by 
$g\mapsto \sigma(g)^n/\det(g)$ is a square-root of
the defining character $\sigma$ of $L_Y$, 
hence $\kappa$ defines a line bundle
$K$ on $Y$ together with an isomorphism
$\varphi : K^2 \to L_Y$. 
Hence the triple $\T \otimes \pi^*(K^{-1})$
takes the form $(V,\O_N,q)$.
In particular, the triple $\T_{N-Y} \otimes \pi^*(K^{-1})$
on $N-Y$ has its structure group reduced from 
$GO(2n)$ to $O(2n)$. Thus in the special case
where $X$ is the total space of a rank $2$ real vector
bundle on $Y$ with $Y\subset X$ the zero section,
when studying mild degenerations in
which the generic rank is even, it is enough to 
consider quadric invariants only in $PH^*(BO(2n))$.

\medskip

\centerline{\bf Examples}

In the following explicit examples, 
the multiplicity $|\nu|_Y(\det(Q))$ is 
assumed to be odd, as otherwise $\pa$ 
maps each invariant to $0$ as shown.

\example 
{\bf Rank $3$ degenerating to rank $2$}
In this case, the ring of quadric invariant is  
$PH^*(BGO(3)) = \FF_2[\wh{w}_2,\wh{w}_3]$, while
$H^*(BGO(2)) = \FF_2[\lambda, a_1, b_4]/(\lambda a_1)$.
Using the Corollary \ref{odd to even} 
we explicitly determine the following.
\begin{eqnarray*} 
\partial((\wh{w}_2^{2i}\,  \wh{w}_3^{2j})  (Q_{X-Y}))& = &0\\
\partial((\wh{w}_2^{2i+1}\,\wh{w}_3^{2j})  (Q_{X-Y}))& = &
         ((a_1^4 + b_4)^i a_1^{2j+1} b_4^j)(\T^{Q_Y}) \\
\partial((\wh{w}_2^{2i}\,  \wh{w}_3^{2j+1})(Q_{X-Y}))& = &
        ((a_1^4 + b_4)^i  a_1^{2j+2} b_4^j)(\T^{Q_Y}) \\
\partial((\wh{w}_2^{2i+1}\,\wh{w}_3^{2j+1})(Q_{X-Y}))& = &
        ((a_1^4 + b_4)^i  a_1^{2j+4} b_4^j)(\T^{Q_Y})
\end{eqnarray*}

\example 
{\bf Rank $4$ degenerating to rank $3$} 
In this case the quadric invariants of $Q_{X-Y}$, by Proposition 4.3,
are the polynomials in 
$\lambda (Q_{X-Y})$, $a_1(Q_{X-Y})$, $(a_1\,a_3+b_4)(Q_{X-Y})$ and 
$(a_3^2+a_1\,d_{\{1,2\}})(Q_{X-Y})$. 
The characteristic classes of the triple $\T^{Q_Y}$ 
are the polynomials in $c(\T^{Q_Y})$, $\widehat{w}_2(\T^{Q_Y})$, and
$\widehat{w}_3(\T^{Q_Y})$. 

Using the Corollary \ref{even to odd} 
we can explicitly determine the following
\begin{eqnarray*} 
\partial(\lambda(Q_{X-Y}))       & = & 0 \\
\partial(a_1(Q_{X-Y}))           & = & 1 \\ 
\partial((a_1\,a_3+b_4)(Q_{X-Y}))& = & \widehat{w}_3 (\T^{Q_Y}) \\
\partial((a_3^2+a_1\,d_{\{1,2\}})(Q_{X-Y}))& = &
        (c\widehat{w}_3 +\widehat{w}_2\widehat{w}_3)(\T^{Q_Y}) 
\end{eqnarray*}

\example 
{\bf Rank $5$ degenerating to rank $4$}
In this case, the quadric invariants of $Q_{X-Y}$ are
polynomials in $\wh{w}_2(Q_{X-Y}),\ldots,\wh{w}_5(Q_{X-Y})$.
The boundary behavior of the $\wh{w}_i(Q_{X-Y})$,
following Corollary \ref{odd to even}, is 
given in Example \ref{wi} below.
As another example, 
$\pa((\wh{w}_2^2\wh{w}_3)(Q_{X-Y})) = (a_1^2 b_4)(\T^{Q_Y})$.  

\example 
{\bf Rank $6$ degenerating to rank $5$} 
In this case the quadric invariants of $Q_{X-Y}$, by Theorem 4.8,
are the polynomials in the cohomology classes
{\small
$$\lambda (Q_{X-Y}),\,{\alpha '}_{1}(Q_{X-Y}),\,
{\alpha'}_{3}(Q_{X-Y}),\,{\alpha '}_{5}(Q_{X-Y}), \,
{\beta '}_{8}(Q_{X-Y}),\,{\beta '}_{12}(Q_{X-Y}),\,
{\delta '}_{\{2,3\}}(Q_{X-Y}).$$
}
The characteristic classes of the triple $\T^{Q_Y}$ 
are the polynomials in $c(\T^{Q_Y})$, $\widehat{w}_2(\T^{Q_Y})$, 
$\widehat{w}_3(\T^{Q_Y})$, $\widehat{w}_4(\T^{Q_Y})$, and 
$\widehat{w}_5(\T^{Q_Y})$. 
Using the Corollary \ref{even to odd}
we get the following.
\begin{eqnarray*} 
\partial(\lambda(Q_{X-Y}))      & = & 0 \\
\partial({\alpha '}_1(Q_{X-Y})) & = & 1 \\ 
\partial({\alpha '}_3(Q_{X-Y})) & = & (c+\widehat{w}_2)(\T^{Q_Y}) \\
\partial({\alpha '}_5(Q_{X-Y})) & = &
            (c\,\widehat{w}_2 +\widehat{w}_2^2)(\T^{Q_Y})\\
\partial({\beta '}_8(Q_{X-Y}))  & = &
      (\widehat{w}_3\,\widehat{w}_4+ c^2 \,\widehat{w}_3+
               \widehat{w}_2^2\,\widehat{w}_3)(\T^{Q_Y}) \\
\partial({\beta '}_{12}(Q_{X-Y}))  & = & 
     (c^2\,\widehat{w}_2\,\widehat{w}_5+
       \widehat{w}_2\,\widehat{w}_4\widehat{w}_5+
       c\,\widehat{w}_2^3 \,\widehat{w}_3+
       \widehat{w}_2^2\,\widehat{w}_3\,\widehat{w}_4)(\T^{Q_Y}) \\
\partial({\delta '}_{\{2,3\}}(Q_{X-Y})) & = & 
 (c^3\,\widehat{w}_2+\widehat{w}_2^4
+\widehat{w}_5\,\widehat{w}_3)(\T^{Q_Y})
\end{eqnarray*}

\example\label{wi} 
We consider the general case of odd degenerating to even, where
the quadric invariants of $Q_{X-Y}$ are all the polynomials 
in $\wh{w}_2(Q_{X-Y}),\ldots, \wh{w}_{2n+1}(Q_{X-Y})$.
The Corollary \ref{odd to even} in particular gives 
\begin{eqnarray*} 
\partial(\wh{w}_{2r}(Q_{X-Y}))& = &
\sum_{i=1}^{r}{{2n+1-2i}\choose{2r-2i}}(a_1^{2r-2i}a_{2i-1})(\T^{Q_Y})\\
\partial(\wh{w}_{2r+1}(Q_{X-Y}))& = &
\sum_{i=1}^{r}{{2n+1-2i}\choose{2r-2i}}(a_1^{2r-2i+1}a_{2i-1})(\T^{Q_Y})\\
\end{eqnarray*}
Note that $\partial(\wh{w}_{2r+1}(Q_{X-Y})) = 
a_1(\T^{Q_Y}) \partial(\wh{w}_{2r}(Q_{X-Y}))$, which 
reflects the well-known relation 
$\sq^1w_{2i} = w_{2i+1} + w_1w_{2i}$ (Wu's formula).



\section{Appendix : Material from Toda [T]}

For any integer $N\ge 1$, let $x_1, \ldots ,x_N$ and $y$ be
independent variables over the field $\FF_2$. Let $A$ be the polynomial
ring $A=\FF_2[x_1,\ldots,x_N]$, and let 
$\phi : A \to \FF_2[t] \otimes A$ be the ring homomorphism
defined by
$\phi(x_r) = \sum_{i=0}^r {N-i \choose r-i}t^{r-i} \otimes x_i$
where $x_0=1$. This makes $A$ into a Hopf algebra comodule over $\FF_2[t]$,
as follows from the topological interpretation of $\phi$ as 
the map $H^*(BG) \to H^*(B\Gamma)\otimes H^*(BG)$ induced by the
multiplication map
$\Gamma \times G \to G$ where $(G,\Gamma)$ is the pair 
$(GL(N,\C), \C\,^*)$ or the pair $(O(N), \{ \pm 1\})$,  
consisting of a topological group $G$ and a central subgroup 
$\Gamma$. Some of the basic properties of $\phi$ 
are just the Hopf algebra comodule axioms. 
The {\bf primitive ring} $PA$ is the subring of $A$ defined by
$ PA = \{x \in A ~|~ \phi(x) = 1\otimes x \}$.
Toda [T] gives a finite set of generators and relations for the
ring $PA$ when $N = 4m +2$ for some $m\ge 0$, and also when $N=4$,
which is recalled below. 

Consider the maps $d_i :A \to A$, for $i \ge 0$, defined   
by $d_i(x_r)={N-r+i \choose i} x_{r-i}$, so that 
$\phi(x) = \sum_{i\ge 0} t^i \otimes d_i(x)$.
In particular, the map $d_1 : A\to A$ is the derivation 
$s = \sum_{1\le i\le N} (N-i+1) x_{i-1} {\pa \over \pa x_i}$
which played an important role in [H-N], with $d_1\circ d_1 =0$. 

The following proposition is an important step in 
Toda's determination of the ring of primitive elements. 

\begin{proposition}\label{main} 
{\rm (Lemma 3.6 of Toda [T]) }
Let $N = qh$ where $q$ is a power of $2$ and $h$ is odd. 
Consider the $\FF_2$-vector subspace (that is, additive subgroup) 
$B$ of $A$, defined by
$B = \{\, a \in A \,|\, d_i(a)=0 \mbox{ for all } i \ge q \,\}$. 
Then the multiplication map $\FF_2[x_q] \times B \to A$ 
which sends $(f(x_q), b) \mapsto f(x_q)b$ induces an 
$\FF_2[x_q]$-linear bijection 
$\FF_2[x_q] \otimes_{\FF_2} B \to A$. 
In particular, it induces an additive isomorphism
$\psi : B \to A/x_qA$. 
\end{proposition}

\rem 
$B$ is not a subring of $A$, that is, the inclusion
$B\hra A$ does not preserve multiplication.
However, the primitive ring $PA$ is a subring of both $B$ and $A$. 
The importance of the ring $B$ is that it helps us get to $PA$.

\rem\label{recursive} 
The proof in [T] of the surjectivity of 
the map $\FF_2[x_q] \otimes B \to A$ is constructive,
that is, given any element in $A$, one can recursively calculate
an explicit preimage in $\FF_2[x_q] \otimes B$. In particular, given
any $a \in A/x_qA$, we can determine explicitly, in a recursive manner,
its preimage $\psi^{-1}a \in B$ under $ \psi : B \to A/x_qA$.

It can be seen that $d_i(x_q) = x_{q-i}$ for all $1\le i\le q-1$.
Consider the elements $\wh{x}_i\in A$ defined by
\begin{eqnarray*}
\wh{x}_i & = & x_i \mbox{ for } 1\le i \le q,\\
\wh{x}_{kq} & = &\psi^{-1}(x_{kq}) \mbox{ for } k\ge 2, \mbox{ and }\\ 
\wh{x}_{kq -i } & = & d_i(\wh{x}_{kq}) \mbox{ for } k\ge 2, \mbox{ and }
1\le i \le q-1.
\end{eqnarray*}

The additive isomorphism of abelian groups
$B \to A/x_qA$ is used by Toda to define a 
ring structure (multiplication operation $*$) 
on the additive group $B$, pulling back the ring structure on
the quotient ring $A/x_qA = \FF_2[x_1,\ldots,
x_{q-1}, x_{q+1},\ldots, x_N]$. Under this, 
$B$ becomes the polynomial ring
$B = \FF_2[\wh{x}_1,\ldots, \wh{x}_{q-1},  \wh{x}_{q+1},\ldots, \wh{x}_N]$
in the $N-1$ variables $\wh{x}_i$ ($i\ne q$) defined above 
(see Proposition 3.7 of [T]).

\example\label{n=2m+1} 
When $N=2n+1$ is odd, the elements 
${\widehat x}_r$ are given by 
${\widehat x}_r =\sum_{i=0}^r {2n+1-i \choose r-i}x_1^{r-i}x_i$.
The primitive ring $PA$, when $N=2n+1$ is odd, is the subring
$\FF_2[{\widehat x}_2\ldots, {\widehat x}_{2n+1}] \subset A$.

The above case of odd $N$ is the case where $q=1$. 
Next we consider the case where $q=2$, that is, 
$N = 2n = 4m +2$ is congruent to $2$ modulo $4$.
The simplest such $N$ is $N=2$, 
where it can be seen that $\wh{x}_1 = x_1$ and $PA = B = \FF_2[x_1]$.

\example\label{s and t recursion} 
When $N = 4m +2$ for some
$m\ge 0$, we have $\wh{x}_1 = x_1$, and the elements 
$\wh{x}_3, \ldots, \wh{x}_N \in B$
can be calculated as follows.
To begin with, we define $s_0=1$ and $t_0 =0$. Now for each $i \ge 0$ 
we define $s_i$ and $t_i$ recursively by $s_i=x_2\,t_{i-1}$ and 
$t_i= s_{i-1}+x_1\,t_{i-1}$. For example, $s_1=0$ and $t_1 =1$, and
$s_2=x_2$ and $t_2 =x_1$. 
The elements $\wh{x}_1,\ldots, {\widehat x}_N$ 
are given by 
$${\widehat x}_{2r} = \sum_{i=0}^{2r} 
{4m+2-2r+i \choose i}\,x_{2r-i}\,s_i \mbox{ and } 
{\widehat x}_{2r-1} = \sum_{i=0}^{2r} 
{4m+2-2r+i \choose i}\,x_{2r-i}\,t_i$$

\bigskip
 
\centerline{{\bf The ring $PA$ when $N = 2n = 4m+2$}}

In this case the action can be written down as follows 
(recall that $\wh{x}_1 = x_1$ and $\wh{x}_2 = x_2$ in this case).
$$
\begin{array}{lll}
\phi(\wh{x}_{2k-1})     & = & 1 \otimes \wh{x}_{2k-1} \mbox{ for all }
                              1\le k \le 2m +1 ,\\
\phi(\wh{x}_2)          & = & 1 \otimes \wh{x}_2 + 
                              t \otimes \wh{x}_1 + t^2 \otimes 1, 
                              \mbox{ and }  \\ 
\phi({\widehat x}_{2k}) & = & 1 \otimes {\widehat x}_{2k}  + 
                              t \otimes {\widehat x}_{2k-1} \mbox{ for all }
                              2 \le k \le 2m + 1. 
\end{array}
$$

We now write Toda's system of
generators $\alpha _{2k-1}$, $\delta_T$, and 
$\beta _{4k}$ for the ring $PA$.

For $1 \le k \le 2m+1$, let  
$$\alpha _{2k-1}= \wh{x}_{2k-1} = d_1(\wh{x}_{2k})$$

The ring multiplication $* : B\times B \to B$ can be extended 
as a binary operation $* : A\times A \to A$ by putting
$$b*c = bc + d_1(b)d_1(c)x_2 \mbox{ for all } b,\,c\in A$$
Using this, more generally for any subset 
$T=\{p_1,\ldots,p_r\} \subset \{ 2,\ldots,2m+1 \}$ 
of cardinality $r >1$, we define
$${\delta}_T= d_1({\widehat x}_{2p_1}* \cdots * {\widehat x}_{2p_r})$$

Next, we define for $2 \le k \le 2m+1$,
$$\beta _{4k}= {\widehat x}_{2k} * {\widehat x}_{2k}
+ x_1 {\widehat x}_{2k-1} {\widehat x}_{2k}  $$

\begin{proposition}\label{n=4m+2} 
{\rm (Proposition 3.11 of Toda [T]) } 
Let $N=4m +2$ for $m\ge 0$. Then the primitive ring $PA$ is 
generated by the elements
$\alpha _{2i-1},\,\beta _{4i},\,\delta_T$
defined above.
\end{proposition}

(Toda also gives the relations, for which see the above reference).

\bigskip

\centerline {{\bf Primitive ring when $N=4$}}
In the case $N = 4$, Toda proves the following.

{\bf Proposition } 
(Toda [T] Proposition 3.12) 
{\it For $N = 4$, the ring $PA$ is the 
polynomial ring $\FF_2[x_1, d_4, d_6]$, where
$d_4=x_2^2 + x_1x_3$ and $d_6=x_3^2 + x_1^2x_4 + x_1x_2x_3$.
} 

\bigskip
\pagebreak

\centerline{\bf The primitive rings for $BO(N)$ and $BGL(N)$}

The above algebraic material applies to the following topological cases.
Consider a pair $(G,\Gamma)$ consisting of a Lie group $G$ and a central
subgroup $\Gamma \subset G$. Let $\Gamma \times G \to G$ 
be the multiplication
map, sending $(\gamma, g) \mapsto \gamma g$. 
This is a group homomorphism
as $\Gamma$ is central, so induces a map 
$B(\Gamma \times G) = B\Gamma \times BG \to BG$ 
on the classifying spaces. 
Let $\phi : H^*(BG) \to
H^*(B\Gamma) \otimes H^*(BG)$ be the induced 
ring homomorphism on cohomology.
Now we consider two cases, as follows.

{\bf Orthogonal group } When $(G,\gamma) = (O(N), \{\pm 1\})$,  
$H^*(BG) = \FF_2[w_1,\ldots, w_N]$ is a 
polynomial ring on the Stiefel-Whitney
classes $w_i$, and $H^*(B\Gamma) = \FF_2[w]$. The map
$\phi$ is given by 
$w_r \mapsto \sum_{0\le i\le r} {N-i \choose r-i} w^{r-i} w_r$.

{\bf General linear group } When $(G,\gamma) = (GL(N), \C\,^*)$,  
$H^*(BG) = \FF_2[\ov{c}_1,\ldots, \ov{c}_N]$ 
is a polynomial ring on the mod $2$ 
Chern classes
classes $\ov{c}_i$, and $H^*(B\Gamma) = \FF_2[t]$. The map
$\phi$ is given by 
$\ov{c}_r \mapsto \sum_{0\le i\le r} {N-i \choose r-i} t^{r-i}\ov{c}_r$.

Hence with appropriate change of notation, the material above gives the 
primitive rings for the orthogonal groups 
or general linear groups in the 
cases when $N$ is odd, or 
when $N = 4$ or when $N = 4m+2$ for some $m\ge 0$.

\medskip

\centerline {\bf Notation for generators of $PH^*(BO(N))$ 
and $PH^*(BGL(N))$} 

For $N=4m+2$, we denote the generators for $PH^*(BO(N))$ by
$\alpha_{2k-1} = \wh{w}_{2k-1}$, ${\delta}_T= 
d_1({\widehat w}_{2p_1}* \cdots * {\widehat w}_{2p_r})$
and $\beta _{4k}= {\widehat w}_{2k} * {\widehat w}_{2k}
+ w_1 {\widehat w}_{2k-1}{\widehat w}_{2k}$. 
On the other hand, for $N=4m+2$ we 
will denote the generators for $PH^*(BGL(N))$ by
$\alpha_{2k-1}'' = \wh{c}_{2k-1}$, ${\delta}_T''= 
d_1({\widehat c}_{2p_1}* \cdots * {\widehat c}_{2p_r})$
and $\beta _{4k}''= {\widehat c}_{2k} * {\widehat c}_{2k}
+ \ov{c}_1 {\widehat c}_{2k-1} {\widehat c}_{2k} $.

\rem 
The ring homomorphism 
$\theta^*: H^*(BGL(N)) \to H^*(BO(N))$ induced by the
inclusion $\theta: O(N) \to GL(N)$ maps 
$\ov{c}_i \mapsto w_i^2$ for each $i$. 
Hence under $\theta^*:H^*(BGL(4m+2)) \to H^*(BO(4m+2))$ 
the invariants ${\alpha''}_{2i-1}$,
${\beta''}_{4i}$ and ${\delta''}_T$ map to ${\alpha }_{2i-1}^2$,
${\beta }_{4i}^2$ and ${\delta }_T^2$ respectively.

\medskip

\centerline{\bf The tensor product 
$\C\,^2 \otimes \C\,^{2m+1} \to \C\,^{4m+2}$}

Consider the homomorphism 
${\tau''}: GL(2m+1)\times GL(2)\to GL(4m+2)$
defined by the tensor product
$\C\,^2 \otimes \C\,^{2m+1} \to \C\,^{4m+2}$. This induces a 
ring homomorphism
$$B(\tau'')^*:H^*(BGL(4m+2))\to H^*(BGL(2))\otimes H^*(BGL(2m+1))$$
The following proposition is a part of the Proposition
4.2 in Toda [T].

\begin{proposition}\label{tau''=0} 
With the above notations we have
$B(\tau'')^*({\alpha''} _{2i-1})=0$ for all $i>1$, 
and $B(\tau'')^*({\delta''} _{T})=0$ for every 
$T=\{i_1,\ldots,i_k\}\subset \{2,\ldots,2m+1\}$.
\end{proposition}

This has the following analog for $O(4m+2)$.

\begin{proposition}\label{tau=0} 
Let the group homomorphism 
$\tau : O(2) \times O(2m+1) \to O(4m+2)$ be induced by the 
tensor product. Then we have
$B(\tau)^*(\alpha _{2i-1})=0$ for all $i>1$, and 
$B(\tau)^*(\delta _{T})=0$ for every 
$T=\{i_1,\ldots,i_k\}\subset \{2,\ldots,2m+1\}$.
\end{proposition}

\proof Observe that the following diagram commutes
$$
\begin{array}{clc}
H^*(BGL(4m+2))&\stackrel{B({\tau''})^*}{\to} & 
H^*(BGL(2))\otimes H^*(BGL(2m+1))\\
{\scriptstyle \theta^*} \downarrow 
& & \downarrow {\scriptstyle \theta^* \otimes \theta^*} \\
H^*(BO(4m+2))&\stackrel{B(\tau)^*}{\to}& H^*(BO(2))\otimes H^*(BO(2m+1))
\end{array}
$$
where $\theta^*: H^*(BGL(4m+2))\to H^*(BO(4m+2))$ is the natural map 
which takes the element $\ov{c}_i \mapsto w_i^2$.
Hence this map also takes the elements 
${\wh c}_{2i-1}$ to ${\wh w}_{2i-1}^2$
and ${\delta''}_{T}$ to $\delta _{T}^2$.
Now for any $a \in H^*(BGL(4m+2))$ such that $B({\tau''})^*(a)=0$, we have 
$B(\tau )^*\theta^*(a)=0$.
This implies that $B(\tau )^*({\wh w}_{2i-1}^2)=0$ 
and $B(\tau )^*(\delta _{T}^2)=0$.
Now the proposition follows from the fact that
$H^*(BO(2))\otimes H^*(BO(2m+1))$ 
has no zero divisors, being a polynomial ring. \hfill $\square$


\section*{References} \addcontentsline{toc}{section}{References} 

[H-N] Holla, Y.I. and Nitsure, N. : Characteristic Classes for $GO(2n)$.
To appear in the Asian Journal of Mathematics.
(e-print math.AT/0003147 available on xxx.lanl.gov).

[N] Nitsure, N. : Topology of Conic Bundles.
J. London Math Soc. (2) {\bf 35} (1987) 18-28. 

[T] Toda, H. : Cohomology of Classifying Spaces. In the volume
{\sl Homotopy Theory and Related Topics} (Editor H. Toda),   
Advances in Pure Math. {\bf 9} (1986) 75-108.

\bigskip 

\bigskip 

\bigskip 

{\it 
Address : School of Mathematics, Tata Institute of Fundamental
Research, Homi Bhabha Road, Mumbai 400 005, India.

e-mails: yogi@math.tifr.res.in, nitsure@math.tifr.res.in 
} 

\bigskip 

\centerline{15-VIII-2000, revised on 11-V-2001} 

\end{document}